\documentclass[reqno,11pt]{amsart}
\usepackage[utf8]{inputenc}
\usepackage{amsmath, latexsym, amsfonts, amssymb, amsthm, amscd, color,caption}
\usepackage[wide,rightcaption]{sidecap}
\usepackage{multirow,stmaryrd}

\usepackage{graphics,epsf,psfrag}
\numberwithin{equation}{section}

\newtheorem{theorem}{Theorem}[section]
\newtheorem{lemma}[theorem]{Lemma}
\newtheorem{proposition}[theorem]{Proposition}

\newtheorem{rem}[theorem]{Remark}

\newcommand{\ind}{\mathbf{1}}

\newcommand{\R}{\mathbb{R}}

\renewcommand{\tilde}{\widetilde}

\newcommand{\cF}{{\ensuremath{\mathcal F}} }

\newcommand{\cN}{{\ensuremath{\mathcal N}} }

\newcommand{\cW}{{\ensuremath{\mathcal W}} }

\DeclareMathSymbol{\leqslant}{\mathalpha}{AMSa}{"36} 
\DeclareMathSymbol{\geqslant}{\mathalpha}{AMSa}{"3E} 
\DeclareMathSymbol{\eset}{\mathalpha}{AMSb}{"3F}     
\renewcommand{\leq}{\;\leqslant\;}                   
\renewcommand{\geq}{\;\geqslant\;}                   
\newcommand{\dd}{\,\text{\rm d}}             

\newcommand{\bbC}{{\ensuremath{\mathbb C}} }

\newcommand{\bbE}{{\ensuremath{\mathbb E}} }

\newcommand{\bbN}{{\ensuremath{\mathbb N}} }

\newcommand{\bbP}{{\ensuremath{\mathbb P}} }

\newcommand{\bbR}{{\ensuremath{\mathbb R}} }
\newcommand{\bbS}{{\ensuremath{\mathbb S}} }

\newcommand{\ga}{\alpha}
\newcommand{\gb}{\beta}
\newcommand{\gd}{\delta}
\newcommand{\gep}{\varepsilon}       

\newcommand{\gP}{\Phi}

\newcommand{\go}{\omega}

\newcommand{\gl}{\lambda}

\newcommand{\gs}{\sigma}

\makeatletter
\def\captionfont@{\footnotesize}
\def\captionheadfont@{\scshape}

\long\def\@makecaption#1#2{%
  \vspace{2mm}
  \setbox\@tempboxa\vbox{\color@setgroup
    \advance\hsize-6pc\noindent
    \captionfont@\captionheadfont@#1\@xp\@ifnotempty\@xp
        {\@cdr#2\@nil}{.\captionfont@\upshape\enspace#2}%
    \unskip\kern-6pc\par
    \global\setbox\@ne\lastbox\color@endgroup}%
  \ifhbox\@ne 
    \setbox\@ne\hbox{\unhbox\@ne\unskip\unskip\unpenalty\unkern}%
  \fi
  \ifdim\wd\@tempboxa=\z@ 
    \setbox\@ne\hbox to\columnwidth{\hss\kern-6pc\box\@ne\hss}%
  \else 
    \setbox\@ne\vbox{\unvbox\@tempboxa\parskip\z@skip
        \noindent\unhbox\@ne\advance\hsize-6pc\par}%
\fi
  \ifnum\@tempcnta<64 
    \addvspace\abovecaptionskip
    \moveright 3pc\box\@ne
  \else 
    \moveright 3pc\box\@ne
    \nobreak
    \vskip\belowcaptionskip
  \fi
\relax
}
\makeatother
\def\writefig#1 #2 #3 {\rlap{\kern #1 truecm
\raise #2 truecm \hbox{#3}}}

\newcommand{\dist}{\mathrm{dist}}

\begin{document}

\title[phase  diffusion in stochastic limit cycle oscillators]{Small noise and long time phase  diffusion in stochastic limit cycle oscillators}

\author{Giambattista Giacomin}
\address{Universit\'e Paris Diderot, Sorbonne Paris Cit\'e,  Laboratoire de Probabilit{\'e}s et Mod\`eles Al\'eatoires, UMR 7599,
            F-75205 Paris, France}

\author{Christophe Poquet}
\address{Universit\'e de Lyon, Universit\'e Lyon 1, Institut Camille Jordan, UMR 5208,
F-69622 Villeurbanne, France}

\author{Assaf Shapira}

\begin{abstract}
We study the effect of additive Brownian noise on an ODE system that has a stable hyperbolic limit cycle, for initial data that are attracted to the limit cycle. The analysis is performed in the limit of small noise -- that is, we modulate the noise by a factor $\gep \searrow 0$ -- and on a long time horizon. We prove explicit estimates
on the proximity of the noisy trajectory and the limit cycle up to times $\exp\left(c \gep^{-2}\right)$, $c>0$, and we show both that on the time scale $\gep^{-2}$ the \emph{dephasing} (i.e., the difference between noiseless and  noisy system measured in a natural coordinate system that involves a phase) is close to a Brownian motion with constant drift, and that on longer time scales the dephasing dynamics is dominated, to leading order, by the drift. The natural choice of coordinates, that reduces the dynamics in a neighborhood of the cycle to a rotation, plays a central role and makes the connection with the applied science literature in which noisy limit cycle dynamics are often reduced to a diffusion  model for the phase of the limit cycle. 
 \\[10pt]
  \\[10pt]
\end{abstract}

\date{\today}

\maketitle

\section{Introduction}

\subsection{Noise induced dephasing phenomena}
Periods, cycles, rhythms are omnipresent and they play a fundamental role. 
And in fact  dynamical models proposed in a variety of fields display (asymptotically) stable periodic behavior, i.e.  (part of the) trajectories 
are attracted by a periodic trajectory. 
Important examples come from
ordinary differential equations (ODE) with stable limit cycles, like  the   ODE systems  for pray-predator dynamics \cite[Ch.~3]{cf:M}, but  \cite{cf:M} contains several examples from life sciences (gene networks, neural systems,$\ldots$). Of course examples come also from physics, chemistry and other sciences
\cite{cf:Bres,cf:ET,cf:K,cf:PRK,cf:W}. 
It is often  the case that  the ODE model is the result of averaging and/or neglecting plenty of details of the original system that is more faithfully modeled by keeping a huge number of degrees of freedom. Introducing noise is therefore a way to go a step closer to reality. It is then natural to think of the noise as \emph{small}, for example when one is considering the dynamics of
{\sl macroscopic} quantities, i.e. averages of quantities of interest over a whole population. But the question then is:
what is the effect of noise on this type of limit cycles? 

This is of course not a novel question and it has been often tackled aiming at {\sl reducing the system to a phase}.
 It is known in fact  that in absence of noise and in  the proximity of the limit cycle 
such  ODE systems can be reduced to the dynamics of a phase: even more, the system can be mapped to constant speed rotation on the unit circle \cite{cf:iso}. It is therefore natural to seek for phase reductions also in the stochastic setting
and a phase reduction for stochastic systems  is proposed for example in \cite{cf:K} and has been employed in a number of contexts, see for example the references in \cite{cf:YA}. But in \cite{cf:YA} it has been pointed out that the stochastic phase reduction  model that has been used is not accurate and that the noise, even when it is white, induces a {\sl frequency shift}.  
In  \cite{cf:YA} a formal small noise development of the solution is given: of course, since the noise is weak the leading order behavior -- what we may call the {\sl macroscopic behavior} -- is just the noiseless behavior. The purpose of   \cite{cf:YA} and of much of the  literature -- similar analyses in fact are developed for example in \cite[Ch.~6]{cf:Bres} and
\cite[\S~10.2]{cf:Scott}, with plenty of references -- has been on catching the next order correction. 
Our purpose is to put these works on rigorous and more  quantitative grounds, changing somewhat the perspective. The question is rather: on which time scale the difference between the phase dynamics in the noisy and noiseless systems becomes macroscopic and, on this time scale, what is the dynamics? 
The answer is that the scale is $\gep^{-2}$ and the dephasing dynamics is a diffusion, in fact a Brownian motion with a constant drift -- the noise induced frequency shift of  \cite{cf:YA}. We will also aim at   longer time scales -- in fact, till the Large Deviations scale in which the noise may induce escapes from the limit cycle  -- and will show that the noise induced frequency shift dominates the phase dynamics. It is worth pointing out that the phase diffusion result we analyze is a macroscopic effect of a microscopic noise that happens on much shorter times scales than the Large Deviations scale, even if in much of the literature macroscopic effects of microscopic noise 
are often identified with Large Deviations phenomena, see   e.g. \cite[Ch.~9]{cf:Scott}.

\subsection{ODEs, limit cycles and noise}
Consider the
Ordinary Differential Equation (ODE)
\begin{equation}
\label{eq:det}
 \dot x_t \, =\, F(x_t) \, ,
\end{equation}
supplied by the initial condition $x_0\in \R^d$. 
In order to have a well defined evolution at least locally in time 
we will assume that $F: \R^d \to \R^d$ is locally Lipschitz.
In reality our main assumption is that \eqref{eq:det} has a hyperbolic (locally) stable 
limit cycle $M:=\{q_t:\, t\in [0, T)\}$: in absence of noise, trajectories are locally attracted to the limit cycle and their distance to the cycle decays exponentially fast, see Section~\ref{sec:model} for the precise set-up.

We then consider
the (strong) solution $X^\gep_\cdot$
 of the Stochastic Differential Equation (SDE)
\begin{equation}
\label{eq:main}
 \dd X^\gep_t\, =\, F\left(X^\gep_t\right)\dd t + \gep \, G\left(X^\gep_t\right) \dd B_t\, ,
\end{equation}
where $\gep>0$,  $G(\cdot)$ is a  locally Lipschitz function from $\R^d$ with values in 
the $d\times m$ matrices and $B_\cdot$ is a standard $m$ dimensional Brownian motion. 
$\Vert A \Vert$ for a $d\times m$ matrix $A$ is the norm of $A$ as linear operator, i.e., with the  choice of the norm
 $\Vert x\Vert:=\max_{j=1, \ldots, d} \vert x_j\vert$ that we keep   throughout the paper,
$\max_i \sum_j \vert A_{i,j}\vert$.
Without loss of generality (see Remark~\ref{rem:Stratonovich}) we interpret $G(X^\gep_t)\dd B_t$ 
in the It\^o sense. These conditions are 
 sufficient to ensure the existence of a unique solution, which may explode in finite time. 
 
It is easy to see for example that if $X^\gep_0=x_0$ sufficiently close to $M$ then for every ${t_f}>0$ and every $\eta>0$ 
\begin{equation}
\label{eq:res-0}
\lim_{\gep \searrow 0}
\bbP 
\left( \sup_{t\in [0, {t_f}]} \left \Vert 
X^\gep_t -x_t \right\Vert  \, > \, \eta \right) \, =\, 0\,.
\end{equation}
In reality this result depends very little on the existence of an hyperbolic  limit cycle, which, for the purpose of \eqref{eq:res-0}, we use  just because it
guarantees that $x_t$ is well defined for every $t>0$: for example, \eqref{eq:res-0} holds simply assuming that $F(\cdot)$ is globally Lipschitz (\cite[Ch.~2]{cf:FW} is dedicated to estimates of this type). A more subtle question is the validity of \eqref{eq:res-0} for longer times, that is if we let ${t_f}$ depend on $\gep$ and become arbitrarily large as $\gep$ becomes small. The fact that the quadratic variations of the (local) martingale  $\gep \, \int_0^{{t_f}}G\left(X^\gep_t\right) \dd B_t$
is $\gep^2 \int_0^{t_f} \left(G\left(X^\gep_t\right)\right)^2 \dd t$ and therefore $O(\gep^2 {t_f})$ (at least if $G$ is bounded),
suggests that the validity of \eqref{eq:res-0} can be pushed up to ${t_f}={t_f}(\gep)=o(\gep^{-2})$ and it certainly suggests
also the breakdown of \eqref{eq:res-0} if  ${t_f} (\gep) =c \gep^{-2}$, any $c>0$. The heuristics of what happens on the time scale $\gep^{-2}$ is just the following: the hyperbolic character of the limit cycle fights against the noisy perturbation
and recalls constantly the trajectory $X^\gep_\cdot$ to $M$. However there is no recalling force in the tangential direction 
to $M$ and the small noise has an effect that becomes macroscopic on times proportional to $\gep^{-2}$ with an effect that 
is qualitatively the same as  the elementary fact that $\gep B_{t\gep^{-2}}$ is again a standard Brownian motion. 
We stress that we used \emph{qualitatively} because a more accurate analysis shows that  in general, i.e. in absence of special symmetries, a drift appears.

We will carry through our analysis and prove a precise and quantitative version of the statement to which we just hinted 
  under the hypothesis that $F(\cdot)$ is $C^2$ in a neighborhood of $M$.
A particular case of  the   results that we are going to prove  
can be stated  by now: there exists a random process $\theta^{\gep}_t$, adapted to the natural filtration 
of $\{B_{t \gep^{-2}}\}_{t\ge 0}$, such that if $X_0^\gep=x_0\in M$ then 
for every $t_f>0$ and $b \in (0,1)$
\begin{equation}
\label{eq:main0}
\lim_{\gep \searrow 0}
\bbP \left( \sup_{t \in [0, t_f]} \left \Vert X_{t \gep^{-2}}^\gep -x_{t\gep^{-2}+ \theta^\gep_t}
\right\Vert > \gep^b \right) \, =\, 0\,, 
\end{equation}
and the law of $\theta^\gep_\cdot$, seen as a random  element of $C^0 ([0, \infty); \bbR)$,
converges as $\gep\searrow 0$ in law \cite[Ch.~2]{cf:Bill} to the Brownian motion with constant drift $\{ \gs B_t+ b t\}_{t\ge 0}$,
where $\gs>0$ and $b \in \R$ depend in a (highly) non trivial way on $F(\cdot)$ and $G(\cdot)$. 
We will give   expressions for  $\gs$ and $b$: for example they can be expressed 
if we know the {\sl Floquet} matrix of the limit cycle we consider, along with $G$, $F$, the Jacobian of $F$  and 
the Hessian of $F_j$, $j=1, \ldots , n$, on the limit cycle. Moreover we will see that a version of \eqref{eq:main0} holds 
also if one requires simply that $x_0$ is close to $M$. 

We will present results also on longer times: \eqref{eq:main0} just deals with times of order $\gep^{-2}$ but if
the Brownian motion with drift $b$, i.e.  $\{ \gs B_t+ b t\}_{t\ge 0}$
is a faithful description of the system on longer time scales then it will be the {\sl noise-induced frequency shift} $b$
(to use the language of \cite{cf:YA}) that dominates. And this is indeed the case, at least until the system
escapes from the domain that is attracted to the limit cycle: this may eventually happen on times
of the order $\exp(c\gep^{-2})$, $c>0$, by effect of Large Deviation events in the driving noise \cite{cf:FW}.

\section{Mathematical set-up and main results}
\subsection{The model and the basics of Floquet Theory}
\label{sec:model}
We restart from \eqref{eq:det} with $F: \R^d \to \R^d$ differentiable -- we prefer $C^1$ to Lipschitz to avoid some minor uninteresting complications -- which is solved by the $T$-periodic 
differentiable function $\bbR \ni t \mapsto q_t\in  \bbR^d$.  
We will use the notation $\Phi(x,t)$ for the solution $x_t$ with $x_0=x$.
As before $M=\{q_t:\, t \in [0, T)\}$, so for example $M=\{\Phi(q_s,t): \, t \ge 0\}$ for every $s$, and we assume that 
$F(\cdot)$ is $C^2$ in a neighborhood of $M$: more precisely we assume that $F(\cdot)$ is $C^2$ in
\begin{equation}
\label{eq:Mdelta}
M_{\gd}\, :=\, \left\{ x:\, \dist(x, M)< \gd\right\}\, ,
\end{equation}
for some $\gd >0$.
 Consider then the evolution 
 linearized  near $M$, that is the linear equation
(called {\sl first variational equation}) 
\begin{equation}
\label{eq:det linearized}
 \dot z_t\, =\, DF(q_t)z_t\, ,
\end{equation}
where $DF(q_t)$ is the Jacobian of $F$ at $q_t$. 
Therefore \eqref{eq:det linearized} is a homogeneous linear equation with periodic coefficients: Floquet theory 
provides a complete understanding of the solutions, even if there are very few cases in which 
\eqref{eq:det linearized} can be solved explicitly. We refer to \cite[Ch.~3]{cf:Teschl} and \cite[\S~2.4]{cf:chicone} for a complete treatment 
of Floquet theory and we give here  some of the features that we need: more can be found in \S~\ref{sec:prelim}.
 
The solution of \eqref{eq:det linearized} starting from $z_s$
at time $s$ is given by
the product $\Pi(t,s)z_s$, where $\Pi(t,s)$ is the {\sl principal matrix solution} associated to the periodic
solution $q_\cdot$, i.e. the solution to
\begin{equation}
\label{eq:eq diff princ mat}
 \partial_t \Pi(t,s)\, =\, D F(q_t)\Pi(t,s)\, ,\quad \Pi(s,s)=I_d\, .
\end{equation}
The principal matrix solution $\Pi(t,s)$ satisfies several basic properties (see \cite[Ch.~3]{cf:Teschl} for
more details): it
is invertible, with inverse $\Pi^{-1}(t,s)=\Pi(s,t)$, and satisfies
the relation
\begin{equation}
 \Pi(t,t_0)\, = \, \Pi(t,s)\Pi(s,t_0)\, .
\end{equation}
Remark in particular that since $z_t= F(q_t)$ is a solution of \eqref{eq:det linearized} starting from $F(q_s)$ at time $s$ we have the relation
\begin{equation}
\label{eq:Feigen}
F(q_t)\, =\, \Pi(t,s) F(q_s)\, ,
\end{equation}
so one is an eigenvalue of $\Pi(t,s)$.
The fact that the solution $q_t$ is periodic of period $T$ implies a periodicity property of the principal matrix solution:
\begin{equation}
\Pi(t+T,s+T)\, =\, \Pi(t,s)\, .
\end{equation}
Other basic properties are in  \S~\ref{sec:prelim}. Here we 
point out that the following representation holds: for every $s$ there exists a matrix $Q(s)$, with $s \mapsto Q(s)$ $T$-periodic, such that
\begin{equation}
\label{eq:later}
\Pi(s+t,s)\, =\, N(s+t,s) e^{-tQ(s)}\, ,
\end{equation} 
$t\mapsto N(s+t,s)$ $T$-periodic and $N(s,s)=I_d$. One then verifies that the matrices
$Q(s_1)$ and $Q(s_2)$ are similar and they have therefore the same eigenvalues and Jordan structure (we will then simply talk of eigenvalues of $Q$). 
Moreover, by \eqref{eq:Feigen}, one sees that $0$ is an eigenvalue. We say that the limit cycle 
$M$ is hyperbolic if only one eigenvalue of $Q$ is equal to $0$, more precisely if its algebraic and geometric multiplicities are equal to one, and we say that $M$ is stable hyperbolic if
it is hyperbolic and if all eigenvalues $\gl_j$ different from $0$ have negative real part: we set 
\begin{equation}
\label{eq:gammaF} 
\gamma_{\textsc f}:= \min_{j:\, \gl_j\neq 0} \Re (-\gl_j)>0\, .
\end{equation}
One can find in  \cite[Ch.~12]{cf:Teschl} the equivalent characterization of (stable) hyperbolic  limit cycles 
in terms of Poincar\'e map and the, also equivalent, interpretation of $M$ as a stable normally hyperbolic manifold.
In particular, if 
$M$ is stable hyperbolic then  if $x \in M_{\gd}$ (the value of $\gd>0$ may be different than the one chosen when we first introduced \eqref{eq:Mdelta})
we have that for every $\gamma< \gamma_F$
\begin{equation}
\label{eq:contract0}
\lim_{t\to \infty}
e^{\gamma t}\dist \left( \Phi(x,t) , M\right)\, =\, 0\, .
\end{equation}

\medskip

\subsection{Isochrons and isochron map}

A much more refined version of \eqref{eq:contract0} holds: for $x\in M$, hence $\Phi(x,t)=q_{t+t_0}$ for some $t_0$, we introduce 
\begin{equation}
\label{eq:W}
W(x)\, :=\,\left\{
y\in \R^d:\, \lim_{t \to \infty}\left \Vert \Phi(y,t)-\Phi(x,t) \right \Vert \, =\, 0 \right \}\, ,
\end{equation}
and, since $M$ is stable hyperbolic, for every $\gamma\in (0,\gamma_F)$
\begin{equation}
\label{eq:W1}
W(x)\, :=\,\left\{
y\in\R^d:\, \sup_{t \ge 0}e^{\gamma t}\left \Vert \Phi(y,t)-\Phi(x,t) \right \Vert \, < \, \infty \right \}\, ,
\end{equation}
We say that $W(x)$ is the {\sl isochron} of $x$. Note that in the standard dynamical systems terminology, $W(x)$ is the stable manifold of $x$ for the map $y\mapsto \Phi(y, T)$.
\smallskip

Here is an important result:
\medskip

\begin{theorem} [Theorem~4.1 in \cite{cf:HPS}] 
\label{th:iso}
For every $x \in M$, $W(x)$ is a $(d-1)$-dimensional  manifold transverse to $F(x)$ at $x$ and of the same regularity as $F(\cdot)$ (hence, $C^k$ if $F(\cdot)$ is $C^k$, $k=1,2, \ldots$). The collection of disjoint open sets $\{W(x)\}_{x \in M}$ is a foliation of the stable manifold $W$  of $M$, that is $W
= \cup_x W(x)$, and $W$ is also open.  Moreover if for every $y\in \R^d$ we call $\theta (y)$ the unique $t$ such that
$y\in W(q_t)$, we have that also $\theta: W \to \bbR \, \mathrm{mod} \,T$ is $C^k$.
\end{theorem}
\medskip

We call $\theta(\cdot)$ \emph {isochron map}: it gives in particular a notion of phase for every point in $W$. It is practical to introduce $\bbS_T:= \bbR \, \mathrm{mod} \,T$
and $\dist_{\bbS_T}(t,s)$ be the arc length between $t$ and $s$. Two important consequences of the definition of isochron map and of Theorem~\ref{th:iso}  are:
\medskip

\begin{itemize}
\item $D\theta (x)F(x)=1$ for every $x\in W$ (here $D\theta(x)$ can be either seen as a differential from acting on $F(x)$
or as the gradient of $\theta(x)$, a row vector,  product the column vector $F(x)$), hence $\frac{\dd}{\dd t} \theta(x_t)=1$ for every $t\ge 0$ and $x_t=\Phi(x_0,t)$
with $x_0 \in W$.
\item Given a compact subset of $W$, there exists a constant $c_{\theta, M}>0$ such that
\begin{equation}
\label{eq:cthetaM}
\left \Vert x - q_{\theta(x)}\right \Vert \, \le \, c_{\theta, M} \text{dist}(x,M)\,,
\end{equation}
for every $x$ in the compact subset we have chosen.
\end{itemize}

\medskip

When we write $r+t$ with $r\in \bbR$ and
$t\in \bbS_T$ we mean $r+ \tilde t$ and $\tilde t$ is the only element of $[0, T)$ such that
$t= \tilde t \, \text{mod}\,T$.
Finally, for every function $f: \R\rightarrow\bbS_T$ we will note by $\tilde{f}$ the lift of $f$, i.e., the unique function $\tilde{f}:\R\rightarrow\R$ such that $\tilde{f}(0)\in[0,T)$, and for all $t$
\begin{equation}
\tilde{f}(t)\mathrm{mod} \,T \,=\,  f(t)\, ,
\end{equation}
see \cite[Prop.~1.33 and Prop.~1.34]{cf:Hatcher}.

\subsection{The stochastic model and main results}
Let  $G: \R^d\to \R^d\times \R^m$ be a Lipschitz continuous (matrix valued) function. By the general theory of SDEs we know how to build a solution to \eqref{eq:main}, but this solution in general explodes in finite time, like its determinist counterpart. Existence and uniqueness up to the (random) explosion time is easily obtained by an approximation and stopping argument (e.g.  
\cite[p.~383]{cf:RY}): one modifies $F(\cdot)$ and $G(\cdot)$ out of the set $\{x \in \bbR:\, \vert x \vert\le n\}$ to make them
globally Lipschitz, then one defines the solution $X^\gep_t$ up to the hitting time of the $\{x \in \bbR:\, \vert x \vert> n\}$ and 
 one finally uses monotonicity of these times in $n$ to pass to the limit. This procedure is useless for us because for the results we are after we can  stop the process upon exiting $W$ or even upon exiting a bounded subset of $W$ (in case $W$ is not bounded), so there is no loss of generality in assuming that $F(\cdot)$ and $G(\cdot)$
are globally Lipschitz (and therefore there exists a unique strong solution $X^\gep_\cdot$ globally in time). Namely 
for every square integrable $\R^d$-valued random variable $X_0^\gep(\cdot)$ measurable with respect to $\cF_0$, $\{\cF_t\}_{t\ge 0}$ is the filtration to which the $n$-dimensional standard Brownian motion $B_\cdot$ is adapted, 
there exists a unique continuous stochastic process $X_\cdot^\gep$   such that
for $j=1, \ldots, d$ and every $t>0$
\begin{equation}
\label{eq:main2}
X^\gep_{j,t}\, =\, X^\gep_{j,0}+ \int_0^t F_j\left(X^\gep_s\right)\dd s + \gep \int_0^t \sum_{l=1}^m
G_{jl}\left(X^\gep_s\right) \dd B_{l,t}\, .
\end{equation}
For $A$ an open subset of $\bbR^d$ we introduce the stopping time $\tau_{\gep, A}:=\inf\{t\ge 0:\, X^{\gep}_t \not \in A \}$. 
 
 \medskip
 
 \begin{rem}
 \label{rem:FW}
 Large Deviation estimates \cite{cf:DZ,cf:FW} can be used to show that 
  if the support of the distribution of $X_0^\gep$ is in $W$ 
  then one can find $c>0$ such that
  $\bbP(\tau_{\gep,W} \ge  \exp(c/\gep^2))$ tends to one as $\gep \searrow 0$. We use this estimate, but for many of the estimates it is useless 
  because while being in $W$ guarantees that the notion of phase makes sense, we need to know that the trajectory is close to 
   $M$ to efficiently approximate the dynamics with a phase dynamics. Of course one could replace $W$ with $M_\gd$,
   $\gd>0$ but small: this improves the situation, but it would induce 
   an error on the phase dynamics that is negligible
   only on times
  $o(1/\gd)$ or, possibly, $o(1/\gd^2)$, that is very far from what we want. For the control of the phase dynamics we are after we need 
  to ensure that the process stays in a \emph{mesoscopic} neighbor of $M$, that is in a neighbor of size that vanishes when $\gep\searrow 0$. These are \emph{moderate deviation} type estimates that we provide explicitly. 
 \end{rem}
 
 \medskip

Our main result is
\medskip

\begin{theorem}
\label{th:main}
Let us choose an arbitrary $t_f>0$ and let us 
assume that there exists $x_0\in W$ such that $\lim_{\gep \searrow 0}X^\gep_0 =x_0$ in probability.
Let us call $\tilde \theta^\gep_t$ the lift $\tilde{\theta\left( X^\gep_{ \cdot\wedge \tau_{\gep, W}} \right)}\left(\gep^{-2} t \right)$. Then for every every $s_\gep= o(\gep^{-2})$ there exists $\eta_\gep=o(1)$
and for every $\gb\in (0,1)$ there exists 
$c_\gb>0$ such that
\begin{equation}
\label{eq:main-prox}
\lim_{\gep \searrow 0}\bbP
\left( \sup_{s\in[0,s_\gep ]}\left \Vert 
X^\gep_{s} -x_s\right \Vert \le \eta_\gep, \, \sup_{t\in [c_\gb \gep^2\vert \log \gep \vert,t_f]} 
\left \Vert
X^\gep_{\gep^{-2}t} -q_{\tilde \theta^\gep_t} \right \Vert \le \gep^\gb
\right)\, =\, 1\, ,
\end{equation}
and the family of processes
\begin{equation}
\label{eq:main-phase-def}
\left\{\tilde \theta^\gep_t
-\theta(x_0)-\gep^{-2}t
\right\}_{t\in [0, t_f]} \, \in C^0\left([0, t_f]; \bbR \right)\, ,
\end{equation}
converges in law as $\gep \searrow 0$ to the process
\begin{equation}
\label{eq:main-phase-lim}
 \big\{ \gs w_t+ bt\big\}_{t \in [0, t_f]}\, ,
\end{equation}
where $w_\cdot$ is a standard Brownian motion
and  the two constants $\gs \ge 0$ and $b \in \bbR$ can be expressed as 
\begin{equation}
\label{eq:sigma1}
\gs^2\, =\, \frac 1T \int_0^T \left( D\theta (q_s) G(q_s)G^{\mathtt{t}}(q_s) D\theta ^{\mathtt{t}} (q_s)   \right) \dd s\, ,
\end{equation}
and
\begin{equation}
\label{eq:b1}
b\, =\, \frac 1{2T} \int_0^T \mathrm{tr}\left( G(q_s) G^{\mathtt{t}}(q_s) D^2\theta(q_s)\right)
\dd s \, ,
\end{equation}
where $D^2\theta$ is the Hessian of $\theta$ and of course \eqref{eq:sigma1} and \eqref{eq:b1} do not depend on 
 $x$.
\end{theorem}

\medskip

A direct consequence of this result is that $\bbP (\tau_{\gep, W}\ge \gep^{-2}t_f)$ tends to one as 
$\gep \searrow 0$, but this is of little interest because, as pointed out in Remark~\ref{rem:FW}, one knows a priori a much stronger result.
As a matter of fact in the proof we use explicit estimates for $\tau_{\gep, M_{\gep^{\gb_0}}}$ 
when dist$(X_0^\gep,M) \le \gep^{\gb_1}$, $0 < \gb_1 < \gb_0$. This type of estimates are in Proposition~\ref{th:prox} and
Proposition~\ref{th:contract} and may be employed to obtain sharper results (or longer time results, like Theorem~\ref{th:main2} below): only very little of their strength is in fact used in our main results -- particularly little in Theorem~\ref{th:main}! -- and Proposition~\ref{th:prox} and
Proposition~\ref{th:contract} can certainly be used as starting point for capturing higher order expansion terms.

\smallskip

It is also worth pointing out that \eqref{eq:main-prox} states two facts: that stochastic and deterministic evolution are almost the same up to times that are $o(\gep^{-2})$ and that it suffices to wait for a time of the order of $\vert \log \gep\vert$ to be 
very close to the limit cycle. The statement can be simplified  if we assume that $X^\gep_0$ is already very close to the limit cycle. Without giving the most general result, we mention that, for example, if $X^\gep_0$ is in $M$ and it is not random \eqref{eq:main-prox} can be replaced  by: for every $\gb\in (0,1)$ 
\begin{equation}
\label{eq:main-prox-M}
\lim_{\gep \searrow 0}\bbP
\left(  \sup_{t\in [0,t_f]} 
\left \Vert
X^\gep_{\gep^{-2}t} -q_{\tilde\theta^\gep_t} \right \Vert \le \gep^\gb
\right)\, =\, 1\, .
\end{equation}

\medskip

\begin{rem}
\label{rem:Stratonovich} 
It is straightforward to generalize Theorem~\ref{th:main}
to the case of
\begin{equation}
\label{eq:main-gen}
 \dd X^\gep_t\, =\, F\left(X^\gep_t\right)\dd t + \gep^2 K\left(X^\gep_t\right)\dd t + \gep \, G\left(X^\gep_t\right) \dd B_t\, ,
\end{equation}
with $K(\cdot)$ a Lipschitz function. The result changes because one has to add to the drift 
the constant term $\frac 1{2T} \int_0^T D\theta(q_s)K(q_s)\dd s$. This observation is relevant because if we were to consider 
the Stratonovich SDE
\begin{equation}
\label{eq:main-strat}
 \dd X^\gep_t\, =\, F\left(X^\gep_t\right)\dd t + \gep \, G\left(X^\gep_t\right) \circ \dd B_t\, ,
\end{equation}
for $G \in C^1$ and $DG$ Lipschitz,  by the standard transformation formula
\eqref{eq:main-strat} can be written in the It\^o form \eqref{eq:main-gen}. The net result is that
Theorem~\ref{th:main} holds with the replacement of \eqref{eq:b1} with
\begin{equation}
\label{eq:b1-strat}
b\, =\, \frac 1{2T} \int_0^T \sum_{i,j,l} G_{i,l}(q_s) D_iD_j\theta(q_s) G_{j,l}(q_s)
  \dd s + \frac 1{2T} \int_0^T \sum_{i,j,l} D_j\theta (q_s) D_i G_{j,l} (q_s) G_{i, l}(q_s) \dd s
 \, ,
\end{equation}
where we have made explicit the formula for the trace for uniformity with the
Stratonovich correction term. Note the use of $D_j$ as partial derivative in the $j^{\textrm{th}}$ coordinate and that $D_i D_j\theta(x)$ is the $(i,j)$ entry of the matrix $D^2 \theta (x)$.
\end{rem}

\begin{rem}
\label{rem:num}
In equations \eqref{eq:sigma1} and \eqref{eq:b1}, $b$ and $\gs$ are given in terms of the derivatives of the isochron map. In fact, these can be expressed explicitly by making a perturbation expansion ear $M$ of the solution to \eqref{eq:det} up to second order. The isochron map $\theta$ is the unique function satisfying
\begin{equation}
\label{eq:theta-oncycle}
\theta(q_t) = t \, \mathrm{mod} \, T,
\end{equation}
and for all $x \in W$
\begin{equation}
\label{eq:theta-oneperiod}
\theta(x) \, =\,  \theta \left( \gP \left( x,T\right) \right)\, .
\end{equation}
If we differentiate \eqref{eq:theta-oncycle} with respect to $t$, recalling that $\dd q_t /\dd t=F\left(q_t\right)$, we  confirm that 
\begin{equation}
\label{eq:theta-first-t-derivative}
D\theta\left({q_t}\right) F(q_t) = 1.
\end{equation}
The second derivative will give a second relation:
\begin{equation}
\label{eq:theta-second-t-derivative}
F^\mathtt{t}\left(q_t\right) D^2\theta\left({q_t}\right) F\left(q_t\right)\,  =\, -D\theta\left({q_t}\right) DF\left({q_t}\right) F\left({q_t}\right)\, .
\end{equation}
Taking a derivative of \eqref{eq:theta-oneperiod} with respect to $x$ at $x=q_t$ on the limit cycle we find a third relation:
\begin{equation}
\label{eq:theta-first-x-derivative}
D\theta\left({q_t}\right) \, =\,  D\theta\left({q_t}\right) D \gP \left( q_t,T \right)\, .
\end{equation}
Finally, differentiating \eqref{eq:theta-oneperiod} twice yields
\begin{equation}
\label{eq:theta-second-x-derivative}
D^2\theta\left({q_t}\right) \,=\,  D\theta\left({q_t}\right) D^2\gP\left(q_t,T\right) + D\gP\left(q_t,T\right)^\mathtt{t} D^2\theta\left({q_t}\right) D\gP\left(q_t,T\right)\, .
\end{equation}
In this notation, $D\theta\left({q_t}\right) D^2\gP\left(q_t,T\right)$ is a matrix whose $i,j$  entry is $\sum_k D_k\theta\left({q_t}\right)D^2_{ij}\gP\left(q_t,T\right)$.
Note that $D \gP(q_t,T)$ is given by the Floquet matrix $e^{-TQ(t)}$, hence it has only one eigenvector, up to normalization, with eigenvalue $1$. By \eqref{eq:theta-first-x-derivative} this eigenvector is exactly $ D\theta\left({q_t}\right) $, and together with the normalization given by \eqref{eq:theta-first-t-derivative} this determines $ D\theta\left({q_t}\right) $ in terms of the Floquet matrix. To find $D^2\theta\left({q_t}\right)$, we  solve the system of linear equations for its coefficients given by \eqref{eq:theta-second-t-derivative} and \eqref{eq:theta-second-x-derivative}. In the Jordan basis of $D\gP\left(q_t,T\right)$ it is straightforward to show that this solution is indeed unique.
\end{rem}

\medskip

Rigorous results in the spirit of Theorem~\ref{th:main}
that have been an important guide toward our result and proof 
can be found in the works about  long time fluctuations of phase boundaries treated 
for Cahn-Allen  SPDEs with bistable symmetric potential in the small noise limit \cite{cf:BdMP,{cf:BBDMP},cf:BBC,cf:F1} and in the zero temperature limit of an interacting Brownian model
\cite{cf:F2}. In  these cases the dynamics takes place close to a hyperbolic invariant manifold, which
  is not a limit cycle: it is a manifold of invariant solutions which is either $\bbR$ or an interval with boundaries (\cite{cf:ABK,cf:BBPS} are also in the same class of problems), and the limit dynamics are diffusions.
 \cite{cf:BGP,cf:Dahms} are about the long time 
dynamics of the mean field plane rotator model: the attracting manifold for the limit PDE is a circle and the resulting dynamics is Brownian motion. All these infinite dimensional results have in common the fact that they
are about \emph{reversible dynamics} (or close to being reversible). This is quite crucial in the analysis and reversibility directly implies the non existence of limit cycles. A  non reversible mean field case, based on the stochastic Kuramoto model, is treated in \cite{cf:LP}: in this case 
the invariant manifold of the limit PDE is still a manifold of stationary solutions, but
the long time dynamics of the system is dominated by a drift
(there is no time scale on which the dynamics is a diffusion). 

With respect to the results we just listed, we require no symmetry properties, neither in the limit cycle, nor in the 
noise (in the sense that $G(\cdot)$ is general). The isochron approach helps in dealing in a rather straightforward way
with the dynamics when the trajectory is close to the limit cycle. To control the proximity of the trajectory to the cycle
-- there is of course a competition between the  drift and the  noise, which is eventually won by noise, unless the cycle is globally stable
 -- we exploit the linear stability theory for periodic trajectories, that is Floquet theory.

\medskip

We now present a result on a  time scale longer than $\gep^{-2}$.
For this note that  as long as $t <  \tau_{\gep,W}$ we can define the winding number $\cW_\gep(t)$ simply 
as the lower integer part of $\tilde \theta^{\gep}/T$.
 In order to avoid stopping let us set $\cW_\gep(t):=\infty$ if $t \ge  \tau_{\gep,W}$. 

\medskip

\begin{theorem}
\label{th:main2}
Under the same assumptions as in Thereom~\ref{th:main} there exists $c>0$ such that
for any choice of $t_\gep>0$ satisfying
\begin{equation}
\label{eq:cond-main2}
\lim_{\gep\searrow 0} \gep^2 t_\gep\, =\, \infty\  \ \ \text{ and } \  \ \ 
\lim_{\gep\searrow 0} \exp(-c\gep^{-2}) t_\gep\, =\, 0\, ,
\end{equation}
we have the following convergence in probability:
\begin{equation}
\label{eq:main2res}
\lim_{\gep \searrow 0 }
\frac{\cW_\gep(t_\gep) -  \left(t_\gep /T\right) }{\gep^2 t_\gep/T }\, =\,  b
\, .
\end{equation}
\end{theorem}
\medskip

One can extract from the proof a value of $c$, but we have chosen not to strive for the optimal value, which 
is of course connected to the solution of a suitable quasipotential problem \cite{cf:FW,cf:DZ}: optimizing and detailing this 
would have
made the arguments longer and heavier. We refer to \cite{cf:Day,cf:BG}, and references therein, for various Large Deviations
issues  related to
 escaping from attracting limit cycles.

 \subsection{Perspectives, a numerical example and organization of the paper}
 Theorem~\ref{th:main2} is much rougher than Theorem~\ref{th:main}. It is natural to conjecture that in order to observe the diffusion effects also on longer time scales, for example on the time scale $\gep^{-\gb}$ with $\gb>2$, one has to subtract 
 a drift term $b_\gep t$, with $\lim_{\gep \searrow 0}b_\gep = b$. In fact our argument of proof allows going somewhat beyond the time scale $\gep^{-2}$, that is we can choose a $\gb$ slightly larger than $2$, because we do prove that
 $b_\gep=b+ O(\gep^a)$ for some $a>0$. But establishing such a result for arbitrary $\gb>2$ necessarily requires the control
 on $b_\gep$ to the adequate order because the drift contribution at time $\gep^{-\gb}$, that  is $b_\gep   \gep^{-\gb}$, has to be be controlled with a better accuracy  than the martingale term which is of the order $\gep^{-\gb/2}$. We believe that this analysis can be performed, but it appears to be extremely cumbersome. In  
\cite{cf:BBDMP} such a result has been proven, but in that case the drift is zero and the diffusive behavior is the leading behavior. 
 
A natural generalization of our work would be to deal with more general stochastic perturbations, for example non Brownian noise: this has been considered in the physical  literature (see e.g.  \cite{cf:Goldobin,cf:Nakao}). And of course our results 
are just in the limit for $\gep \searrow 0$: new phenomena may and do arise for non vanishing noise intensity, see for example
\cite{cf:Newby}.

\begin{SCfigure}[50]
\centering
\leavevmode
\epsfxsize =7.5 cm
\psfragscanon
\psfrag{0}[c][l]{\tiny $0$}
\epsfbox{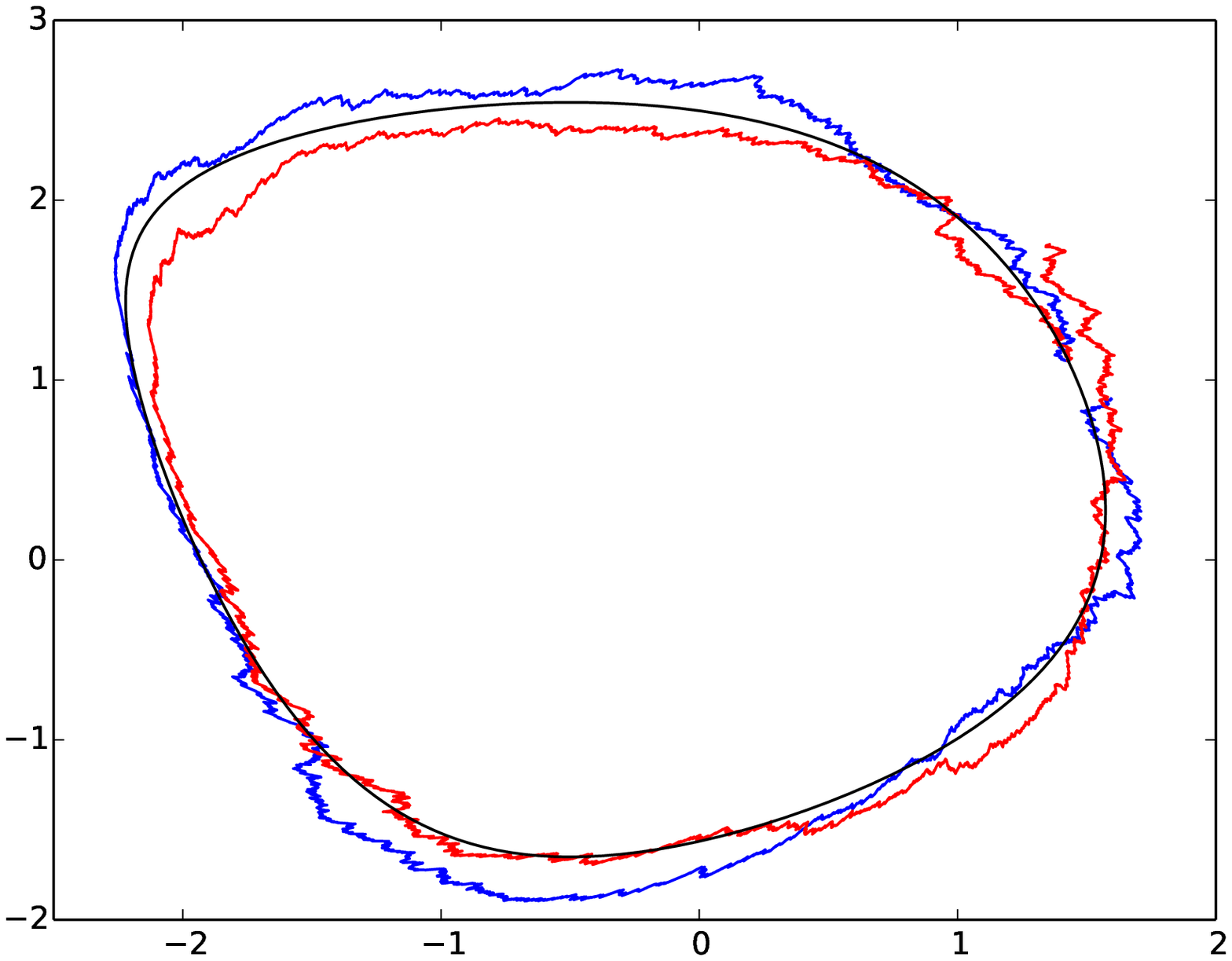}
\caption[.]{\label{fig:three} We plot the limit cycle of a  FitzHugh-Nagumo  system and two trajectories of a stochastic version.  With reference to \eqref{eq:main}: $F_1(x,y)=x-(x^3/3)-y$, $F_2(x,y)=x+1/2$, $\gep=0.1$ and 
\leavevmode\\
\begin{minipage}{\linewidth}
$$
\hspace*{-2.5 cm} 
G(x,y)=
 \frac1{\sqrt{x^2+y^2}}\left(
 \begin{array}{cc}
  x & -y \\
    y & x \\
     \end{array} \right)\left(\begin{array}{cc}
  -1 & -1 \\
    1 & 1 \\
     \end{array}
\right).
$$
  \end{minipage}
The initial condition is on the limit cycle in all cases and the simulations run for one period $T=7.067\ldots$
The numerical data in Table~\ref{table1} refer to this SDE system, with various values of $\gep$.
\phantom{add add add add add add add add add add add add add add add add add add add add add add add add add}
\phantom{add add add add add add add add add add add add add add add add add add add add add add add add add}
}
\end{SCfigure}

To complete our presentation we report some numerical results in Fig.~\ref{fig:three} and Table~\ref{table1}. The numerical evaluation of $\gs$ and $b$
in the $\gep \searrow 0$ limit, that is \eqref{eq:sigma1} and \eqref{eq:b1}, is based on a second order development of the isochrons near $M$: we will not go into the details of the numerical construction of the isochrons that is an issue in its own, 
see e.g. \cite{cf:HdlL} and references therein.

\begin{SCtable}[1.5]
  \begin{tabular}{ | l | l | l | r |}
    \hline
    $\gep$ &  $\gs_{N}$ & $b_N$ & $t_{\text{obs}}$ \\ \hline
    $0.50$ &$1.34 $  &$0.133\pm0.003 $  &$2\, T $\\ \hline
    $0.20$ &$1.23 $  &$0.719\pm0.003 $  &$10\,T $\\ \hline
    $0.10$ &  $1.13$&  $0.699 \pm 0.003$&$40\, T$\\ \hline
    $0.05$&  $1.10$&  $0.690 \pm 0.003$& $160\, T$\\ \hline
    $0.02$ & $1.10$ & $0.689\pm 0.003$ & $1000\, T$\\ \hline \hline
    $0^+$  &  $1.07\ldots$ & $0.688\ldots$ & \text{NA}  \\
    \hline
  \end{tabular}
\caption[.]{\label{table1}
Simulation of the FitzHugh-Nagumo system given in the caption of Figure \ref{fig:three}. 
We have szmpled   $N=50000$ realizations  of the phase at time $t_{\text{obs}}$ and computed the empirical mean $b_N$ and the sample empirical standard deviation $\gs_N$. 
The time $t_{\text{obs}}$ is chosen to keep
$\gep^2 t_{\text{obs}}$ approximately constant, aiming at a comparable error in 
$b_N$: the error is one empirical standard deviation of the mean of the sample.   
In the last line we give $\gs$ and $b$, obtained by evaluating numerically \eqref{eq:sigma1} and \eqref{eq:b1}, see Remark~\ref{rem:num}.
}
\end{SCtable}

Here is the organization of the remainder of the paper:
the proof of Theorem~\ref{th:main} in Section~\ref{sec:proofs}, but it relies 
Proposition~\ref{th:prox} and Proposition~\ref{th:contract} that are in Section~\ref{sec:prox}.
The proof of Theorem~\ref{th:main2} is split in two parts: for \emph{moderately} long times
the proof is in Section~\ref{sec:proofs} and it is then completed in Section~\ref{sec:complete} for times
up to the ones in which Large Deviations events take place.

\section{Proof of the main results}
\label{sec:proofs}

\subsection{Proof of Therem~\ref{th:main}}
The proof has two main steps: early stages and 
long time analysis. 
\subsubsection{Step 1: early stages (approaching $M$)}
Without loss of generality we choose $X_0^\gep=x_0$ non random and in $W$. 
Let us remark from the start that since $x_t=\Phi(x_0, t)$ is attracted by $M$ and since $W$ is open 
there exists $\gd_0>0$, that depends on $x_0$, such that
$\inf_{t\ge 0}\dist (x_t, W^\complement) \ge \gd_0$, and
$x_t\in M_R$ for all $t\ge 0$, $R>0$ depending on $x_0$.
We call $L_F$ the Lipschitz constant of $F(\cdot)$ in the (bounded) set
$M_{R+1}\cap W_{\gd_0/2,-}$, where $W_{\gd_0/2,-}=\{x \in W:\, 
\dist (x, W^\complement) \ge \gd_0/2\}$. We set for this subsection $\tau_\gep:= \tau_{\gep,M_{R+1}\cap W_{\gd_0/2,-}}$ 
so that for every $t\ge 0$

\begin{equation}
\label{eq:fGron1}
\left \Vert  X_{t\wedge{\tau_\gep} }^\gep - x_{t\wedge {\tau_\gep}} \right\Vert \, \le \, L_F \int_0^t 
\left \Vert  X_{s\wedge{\tau_\gep} }^\gep - x_{s\wedge {\tau_\gep}} \right\Vert \dd s + \gep \left \Vert \int_0^{t\wedge {\tau_\gep}} G\left ( X^\gep_s \right) \dd B_s \right \Vert \, ,
\end{equation}
and since, by the Doob's submartingale inequality for every $\eta>0$ and every $t_0>0$
\begin{equation}
\bbP\left( \sup_{t\in [0, t_0]}\gep \left \Vert \int_0^{t\wedge {\tau_\gep}} G\left ( X^\gep_s \right) \dd B_s \right \Vert \, \ge \, \eta\right)\, \le \, 
\frac{t_0 \gep^2}{\eta^2} \Vert G \Vert ^2_{\infty}\, ,
\end{equation}
where $\Vert G \Vert ^2_{\infty}=\Vert G \Vert ^2_{\infty, V_{R+1}\cap W_{\gd_0/2,-}}$,
by Gronwall inequality, from  \eqref{eq:fGron1} we obtain that 
\begin{equation}
\bbP\left( 
\sup_{t \in [0, t_0]} \left \Vert  X_{t\wedge{\tau_\gep} }^\gep - x_{t\wedge {\tau_\gep}} \right\Vert 
  \exp\left(- L_Ft\right)\, \ge\, 
\eta \right)\, \le \,   \frac{t_0 \gep^2}{\eta^2} \Vert G \Vert ^2_{\infty}\, .
\end{equation}
Therefore by choosing $\eta= \gep \vert \log \gep \vert ^b$, $b>1/2$, and $t_0=c_1 \vert \log \gep \vert$ we have that for every $c_1>0$
\begin{equation}
\lim_{\gep \searrow 0}
\bbP\left( \sup_{t \in [0, c_1 \vert \log \gep \vert ]}\left \Vert  X_{t\wedge{\tau_\gep} }^\gep - x_{t\wedge {\tau_\gep}} \right\Vert \, \le \, 
\gep^{1- L_F c_1} \vert \log \gep \vert ^b\right) \, =\, 1\, ,
\end{equation}
and we can therefore remark that since $x_t\in M_R\cap W_{\gd_0,-}$ for every $t$ and since if $c_1<1/L_F$ (which we assume henceforth) then
$X_{t\wedge{\tau_\gep} }^\gep - x_{t\wedge {\tau_\gep}}$ vanishes in probability as $\gep \searrow 0$ we conclude that
$X^\gep_t\in M_{R+1}\cap W_{\gd_0/2,-}$ for every $t \le c_1\vert \log \gep\vert$ with probability approaching one. Therefore  for every $\gb\in (0, 1- L_F c_1)$
\begin{equation}
\label{eq:step1.1}
\lim_{\gep \searrow 0}
\bbP\left( \sup_{t \in [0, c_1 \vert \log \gep \vert ]}\left \Vert  X_{t }^\gep - x_{t} \right\Vert \, \le \, 
\gep^{\gb} \right) \, =\,1\, ,
\end{equation}
which directly implies that  on such an event we  have also
\begin{equation}
\label{eq:step1.2}
 \sup_{t \in [0, c_1 \vert \log \gep \vert] }
\dist_{\bbS_T}\left(  \theta\left(X_{t }^\gep\right) , \theta\left(x_{t}\right) \right) \, \le \, 
L_\theta \gep^{\gb} \, ,
\end{equation}
where $L_\theta$ is the Lipschitz constant of $\theta$ in $V_{R+1}\cap W_{\gd_0/2,-}$.
So \eqref{eq:step1.1}-\eqref{eq:step1.2} directly entail
\begin{equation}
\label{eq:step1.2.1}
\lim_{\gep \searrow 0}
\bbP\left( \sup_{t \in [0, c_1 \vert \log \gep \vert \wedge {\tau_\gep} ]}\left \vert  \tilde{\theta\left(X_{\cdot }^\gep\right)}\left( t\right) - \tilde{\theta\left(x_{\cdot}\right)}\left( t\right) \right\vert \, \le \, 
L_\theta \gep^{\gb} \right) \, =\, 1\, ,
\end{equation}
 Finally, by remarking that $\dist(x_t, M)\le \exp(-\gamma t)$ for $t$ sufficiently large,
$\gamma \in (0, \gamma_{\textsc f})$ (see \eqref{eq:gammaF}-\eqref{eq:contract0}), we have that for $\gep$  small
$\dist(x_{c_1 \vert \log \gep \vert}, M)\le \gep^{\gamma c_1}$ and therefore with probability going to one
\begin{equation}
\label{eq:step1.3}
\dist\left( X^\gep_{c_1  \vert \log \gep \vert}, M\right) \, \le \, \max \left( \gep^\gb, \gep^{\gamma c_1}\right)\, = \, 
\gep ^{\gb_0}\, , 
\end{equation}
where $\gb_0:=\min(\gb, \gamma c_1)$. One can now optimize the choice of $c_1$, by choosing it arbitrarily close to $1/(\gamma+L_F)$, leading to $\gb_0$ smaller but close to $\gamma/(\gamma+L_F)$. In the end, recalling that $\gamma \in (0, \gamma_{\textsc f})$,
$\gb_0$ can be chosen in  $(0,\gamma_{\textsc f}/(\gamma_{\textsc f}+L_F))$.

\subsubsection{Step 2: longer time analysis}
We  can restart the evolution from  $t= c_1  \vert \log \gep \vert$, when, with probability  approaching one,
the trajectory of the process is at a vanishing distance from $M$. 
Therefore from now on  we work with 
\begin{equation}
Y_t^\gep\, := \, X_{t+c_1  \vert \log \gep \vert}^\gep\, .
\end{equation}
 Accordingly, $x_t$ from now is $\Phi (Y_0^\gep, t)$.  
Observe that 
Proposition~\ref{th:prox}  guarantees that $\dist(Y_t^\gep, M)\le \gep ^{\gb_1}$, any $\gb_1<\gb_0$,
with probability approaching one and uniformly on a very long time horizon (well beyond what is sufficient for
the present argument which is just $O(\gep^{-2})$).

Up to now we have controlled the phase of the solution, i.e. \eqref{eq:step1.2}, just in terms of its proximity between 
deterministic and stochastic evolution, i.e. \eqref{eq:step1.1}. On longer time scales  the phase of stochastic 
and deterministic evolutions start differing in a substantial way. We therefore start the sharper analysis of the phase  by 
applying the It\^o formula  to the process
$\tilde{\theta(Y^\gep_{\cdot\wedge \tau_\gep})}-\tilde{\theta(x_{\cdot\wedge \tau_\gep})}$, where we have redefined  $\tau_\gep:= \tau_{\gep ,M_{\gep^{\gb_1}}}$ (this convention will be kept till the end of the proof), obtaining
that 
\begin{multline}
\label{eq:ito1}
u^\gep (t)\, :=\, 
\tilde{\theta \left(Y^\gep_{\cdot\wedge \tau_\gep}\right)}\left( t \right)- \tilde{\theta \left( x_{\cdot\wedge \tau_\gep}\right)}\left( t\right) \, =\\ 
\gep \int_0^t \ind_{[0, \tau_\gep]}(s)D \theta \left(Y^\gep_s\right) 
G\left(Y^\gep_s\right) \dd B_s+
\frac 12 \gep^2 \int_0^t \ind_{[0, \tau_\gep]}(s)\mathrm{tr} \left(
G^{\mathtt{t}}\left(Y^\gep_s\right) D^2 \theta \left(Y^\gep_s\right) G\left(Y^\gep_s\right) 
\right)\dd s
\\
=: \gep \int_0^{t\wedge \tau_\gep} H_1\left(Y_s^\gep\right) \dd B_s+ 
 \gep^2 \int_0^{t\wedge \tau_\gep} H_2\left(Y^\gep_s\right) 
 \dd s\, ,
\end{multline}
where we have used $D\theta F =1 $, $H_1(\cdot)$ is Lipschitz and $H_2(\cdot)$ is $C^0$. Since $Y^\gep_\cdot$ is bound to a neighborhood of 
$M$, we can consider $H_1(\cdot)$  bounded and  Lipschitz with constant
$L_{H_1}>0$, and $H_2(\cdot)$  bounded and uniformly continuous (the continuity modulus will be denoted by
$\go_{H_2}(\cdot)$). 
\subsubsection{Step2.1: replacing $Y^\gep_t$ with $q_{\theta(Y^\gep_t)}$}
We see that for the quadratic variation of the first term in the last line of \eqref{eq:ito1} we have
\begin{multline}
\sup_{t \le \gep^{-2}t_f}\left \vert 
\gep^2 \int_0^{t\wedge \tau_\gep} \left(H_1\left(Y_s^\gep\right)\right)^2 \dd s -
\gep^2 \int_0^{t\wedge \tau_\gep} \left(H_1\left(q_{ \theta (Y_s^\gep)}\right)\right)^2 \dd s
\right \vert \, \le 
\\
 2t_f \Vert H_1 \Vert L_{H_1} \sup_{s \le \gep^{-2}t_f \wedge \tau_\gep} \left \vert Y_s^\gep- q_{ \theta (Y_s^\gep)}\right\vert
\, \le\,  2t_f \Vert H_1 \Vert L_{H_1}c_{\theta, M} \sup_{s \le \gep^{-2}t_f}\textrm{dist}(Y_s^\gep, M)\,,
\end{multline}
where $c_{\theta, M}$ has been introduced  in \eqref{eq:cthetaM}. By Proposition~\ref{th:prox} we see that
this expression is $O(\gep^{\gb_1})$ in probability. Similarly
\begin{equation}
\sup_{t \le \gep^{-2}t_f}\left \vert 
 \gep^2 \int_0^{t\wedge \tau_\gep} \left(H_2\left(Y^\gep_s\right) 
 - H_2\left(q_{ \theta (Y_s^\gep)}\right) \right)
 \dd s\right \vert \, \le \, t_f \go_{H_2}\left(  c_{\theta, M} \sup_{s \le \gep^{-2}t_f}\textrm{dist}(Y_s^\gep, M)\right)\,,
\end{equation}
which, by Proposition~\ref{th:prox}, tends to zero in probability. Of course Proposition~\ref{th:prox} impies also  $\lim_{\gep \searrow 0}\bbP( \tau_\gep> \gep^{-2})=1$.

\subsubsection{Step 2.2: deterministic approximation of $\theta(Y^\gep_\cdot)$ and convergence of drift and quadratic variation}
So to control both the quadratic variation of the first term 
and the second term in the last line of\eqref{eq:ito1} we aim at showing that,  
for every $H$ uniformly continuous,  $\gep^2 \int_0^{t\gep^{-2}\wedge \tau_\gep} H\left(q_{ \theta (Y_s^\gep)}\right) 
 \dd s$ converges, uniformly in $t \in [0, t_f]$, in probability to $t$ times the constant  $\frac 1T \int_0^T H(q_s) \dd s$.
 This is a consequence of the following
 claim: for any $c>0$ 
\begin{equation}
\label{eq:fDoobexp}
\lim_{\gep \searrow 0}\bbP \left( \sup_{n= 1, \ldots, \lfloor  \gep^{-c} \rfloor} \sup_{t \in [(n-1)T,nT]}
 \vert u^\gep(t)-u^\gep((n-1)T)\vert \, > \, \gep \vert \log \gep \vert
\right)\, =\, 0\, .
\end{equation}
In fact, on the complement of the event whose probability is estimated in \eqref{eq:fDoobexp} (choose a $c>2$)
we have 
\begin{equation}
\label{eq:tilde2.2}
\left \vert \tilde{\theta \left(Y^\gep_{\cdot\wedge \tau_\gep} \right)} \left( t \right)- \tilde{\theta \left( x_{\cdot\wedge \tau_\gep}\right)}\left( t\right)
-\left( \tilde{\theta \left(Y^\gep_{\cdot\wedge \tau_\gep}\right)}\left( (n-1)T \right)-
\tilde{\theta \left( x_{\cdot\wedge \tau_\gep}\right) } \left( (n-1)T\right)
\right)
\right \vert\, \le \, \gep \vert \log \gep \vert \, ,
\end{equation}
for every $n\le \gep^{-2}t_f/T$ and every $t \in [(n-1)T, nT]$
and the (random) term between parentheses in the left-hand side of \eqref{eq:tilde2.2} -- we call it $a_{\gep, n}$ -- does not depend on $t$. 
So
we see that on the same event on which 
\eqref{eq:tilde2.2} holds we have also for the same values of $n$ and $t$
\begin{equation}
\label{eq:tilde2.2.1}
\left \vert H \left(q_{ \tilde{\theta \left( Y^\gep_{\cdot\wedge \tau_\gep}\right)} \left( t \right)} \right)- 
H\left(q_{\tilde{\theta \left( x_{\cdot\wedge \tau_\gep}\right)} (t)+a_{\gep,n}}\right)\right \vert \, \le \, \go_H(\gep \vert \log \gep \vert)\, .
\end{equation}
We now observe that 
\begin{equation}
\gep^2 \int_0^{t\gep^{-2}\wedge \tau_\gep}H\left(q_{\theta(Y^\gep_s)}\right) \dd s
\, =\, \sum_{n=1}^{\lfloor (t\gep^{-2}\wedge \tau_\gep)/T\rfloor} \gep^2 \int_{(n-1)T}^{nT} H\left(q_{\theta(Y^\gep_s)}\right) \dd s
+ O(\gep^2)\, ,
\end{equation} 
and using \eqref{eq:tilde2.2.1} (recall also   that $H(q_{\tilde {\theta (y_\cdot)} (t)})= H(q_{\theta(y_t)})$ 
for every $y_\cdot\in C^0([0, \infty); W)$)
we see that 
for $n \le \tau_\gep/T$
\begin{equation}
\label{eq:fr18}
\left \vert \int_{(n-1)T}^{nT} H\left(q_{ \theta (Y_s^\gep)}\right) 
 \dd s -\int_0^T H(\Phi (x,s))\dd s \right \vert \, \le \, T \go_H \left( \gep \vert \log \gep \vert \right)\, , 
\end{equation}
for any $x \in M$.
Therefore, once \eqref{eq:fDoobexp} is established,  we have the convergence we claimed for
$\gep^2 \int_0^{t\gep^{-2}\wedge \tau_\gep} H\left(q_{ \theta (Y_s^\gep)}\right) 
 \dd s$ at the beginning of Step 2.2.

Let us then go back to  \eqref{eq:fDoobexp}
and start by observing that $u^\gep(t) -u^\gep((n-1)T)$ is 
\begin{equation}
\gep \int_{(n-1)T}^t
 \ind_{[0, \tau_\gep]}(s)
 H_1\left(Y^\gep_s\right)
 \dd B_s
+\frac {\gep^2}2  \int_{(n-1)T}^t \ind_{[0, \tau_\gep]}(s)
H_2\left(Y^\gep_s\right)
\dd s\, ,
\end{equation}
and the second term is $O(\gep^2)$ uniformly in $t\in [(n-1)T, nT]$ and in $n$, because the integrand is bounded.
The first term is instead of the form $\gep \int_{(n-1)T}^t Z_s^\gep  \dd B_s$ and $Z^\gep_\cdot$ is a bounded process, so
\begin{equation}
\left\{\exp\left( \gep a \int_{(n-1)T}^t Z_s^\gep  \dd B_s - \frac 12 a^2 \gep^2 \int_{(n-1)T}^t (Z_s^\gep)^2  \dd s\right)\right\}_{t\ge (n-1)T}\, \, 
\end{equation} 
 is a martingale (for every $a\in \bbR$)
and Doob's inequality tells us that 
\begin{multline}
\label{eq:expmartGb1}
\bbE\left[ \sup_{t\in [(n-1)T, nT]} \exp\left( \gep a \int_{(n-1)T}^t Z_s^\gep  \dd B_s\right)\right]
\, \le \\ \exp\left( C_{T,F}\frac{a^2 \gep^2}2
\right) \bbE\left[ \exp\left( \gep a \int_{(n-1)T}^{nT} Z_s^\gep  \dd B_s - \frac 12 a^2 \gep^2 \int_{(n-1)T}^{nT} (Z_s^\gep)^2  \dd s\right)\right] 
\\
=\, \exp\left( C_{T,F}\frac{a^2 \gep^2}2\right)\, ,
\end{multline} 
where $C_{T,F}:=T\Vert H_1 \Vert_\infty^2
$. So by applying the Markov inequality with the choice
$a=\pm \vert \log \gep\vert/(C_{T,F}\gep)$ we obtain
\begin{equation}
\label{eq:expmartGb2}
\begin{split}
\bbP\left( \sup_{t\in [(n-1)T, nT]} \gep \left \vert \int_{(n-1)T}^t Z_s^\gep  \dd B_s\right\vert\, >\, \gep \vert \log \gep\vert\right)
\, &\le \, 2 \exp\left( C_{T,F}\frac{a^2 \gep^2}2- a \gep \vert \log \gep \vert\right)
\\
&=\, 2\exp
\left(-\frac {\vert \log \gep \vert ^2}{2C_{T,F}}
\right)\, ,
\end{split}
\end{equation}
and \eqref{eq:fDoobexp} follows by a union bound estimate.

\subsubsection{Step 2.3: convergence of the process and completion of the proof} 
Since the both the quadratic variation and the drift of
the process 
$t \mapsto u^\gep({\gep^{-2}t})$, cf. \eqref{eq:ito1},
in $C^0([0, t_f]; \bbR)$ converges in probability to the deterministic limits
that we have explicitly identified -- in particular the quadratic variation is linear in time --
the convergence  of the phase process in \eqref{eq:main-phase-def} to the limit Gaussian process 
\eqref{eq:main-phase-def}
 follows by standard arguments (see for example \cite[Ch.~2]{cf:Bill}). The early stages of the evolution do not contribute to the limit 
because the change in phase is $o(1)$.

For what concerns \eqref{eq:main-prox} instead, it is  a matter of exploiting Proposition~\ref{th:prox}
and Proposition~\ref{th:contract}.
More precisely \eqref{eq:main-prox} is in part proven by \eqref{eq:step1.1}, up to times $c_1 \vert \log \gep\vert$
and $\eta_\gep\le \gep^\gb$, $\gb$ identified right before \eqref{eq:step1.1}. 
Starting from such a time, Proposition~\ref{th:prox} that guarantees persistence of proximity to $M$
for times that are even much longer than the ones that we consider here. In particular by the convergence of the
phase process \eqref{eq:main-phase-def} we know that the phase $\tilde \theta^\gep_t$ is close to the phase of the
deterministic solution up to times $t_\gep=o(\gep^{-2})$: this shows that the probability of the first event  in \eqref{eq:main-prox}
tends to one.
 But having gone far enough on the 
$\vert \log \gep\vert$ times scale we can also invoque Proposition~\ref{th:contract} to establish the stronger proximity claimed in the 
second event in \eqref{eq:main-prox}. Therefore the proof of Theorem~\ref{th:main} is complete.
\qed

 \bigskip 
 
The proof of Theorem~\ref{th:main2} is split into two cases: namely we prove it first for $t_\gep$ such that
\begin{equation}
\label{eq:cond-main2.1}
\lim_{\gep\searrow 0} \gep^2 t_\gep\, =\, \infty\  \ \ \text{ and } \  \ \ 
\lim_{\gep\searrow 0} \exp(-\gep^{-\zeta}) t_\gep\, =\, 0\, ,
\end{equation}
with $\zeta$ arbitrarily chosen in $(0,2)$. Then, in Section~\ref{sec:complete}, we will prove it for
\begin{equation}
 \label{eq:cond-main2.2}
\liminf_{\gep\searrow 0} \exp(-\gep^{-\zeta}) t_\gep\, >\, 0\,\  \ \ \text{ and } \  \ \ 
\lim_{\gep\searrow 0} \exp(-c\gep^{-2}) t_\gep\, =\, 0\, .
\end{equation}
\medskip
 
 \noindent
 {\it Proof of Theorem~\ref{th:main2} assuming \eqref{eq:cond-main2.1}}.  
 The early stage dynamics does not have any effect on the result we are after
so we do start with $X^\gep_t$ in $M_{\gep^{\gb_0}}$ for some $\gb_0\in (0,1-\frac{\zeta}{2})$
and therefore, for any $\gb_1\in(0, \gb_0)$, $X^\gep _t \in M_{\gep^{\gb_1}}$
for $t\in [0, \exp(-\gep^{-\zeta})]$ with probability close to one, cf. Proposition~\ref{th:prox}. We therefore restart from
 \eqref{eq:ito1}, with $X^\gep_\cdot$ in place of $Y^\gep_\cdot$,
and we aim at a law of large numbers, in probability, for the process $t \mapsto\tilde{\theta \left(X^\gep_{\cdot\wedge \tau_\gep}\right)}\left( t \right)$.
More precisely we aim at showing that  $(\tilde{\theta \left(X^\gep_{\cdot\wedge \tau_\gep}\right)}\left( t_\gep \right) -t_\gep)/(\gep^2 t_\gep)$
converges to $b$. But $ \tilde{\theta \left( x_{\cdot\wedge \tau_\gep}\right)}\left( t\right)=t+\theta(x_0)$ for every $t \le \tau_\gep$ so we are
back to considering
\begin{equation}\label{eq:ut again}
u_t^\gep\, =\, 
\tilde{\theta \left(X^\gep_{\cdot\wedge \tau_\gep}\right)}\left( t \right)- \tilde{\theta \left( x_{\cdot\wedge \tau_\gep}\right)}\left( t\right)
\end{equation}
with $x_t=\Phi(X^\gep_0,t)$ and the corresponding It\^o formula. The stochastic integral part in the  It\^o formula
gives no contribution to the final result because its quadratic variation is $O(\gep^2 t)$, hence giving a contribution to the
final result that is $O(1/(\sqrt{\gep^2 t_\gep}))=o(1)$, in probability. The drift part in the It\^o formula is treated like in the previous proof after having modified
 \eqref{eq:fDoobexp} by replacing $\gep \vert \log \gep \vert$ with $\gep^{1-b}$, $b\in(\frac{\zeta}{2},1)$ to get estimates on 
 $\exp(\gep^{-\zeta})$ times. The (partial) proof is therefore complete.
\qed

\section{Persistence of proximity and approach to $M$}
\label{sec:prox}

\subsection{The statements}
\label{sec:stat}
In the deterministic case it is not difficult to see that there exists $c_F\ge 1$ such that for $\gd$ sufficiently small, a trajectory that starts in 
$M_\gd$ will not leave $M_{c_F\gd}$. 
Of course this is no longer true in the stochastic setting, but:

\medskip
\begin{proposition}
\label{th:prox}
If $X_0$ is such that $\dist(X^\gep_0, M) \le \gep ^{\gb_0}$ for a $\gb_0\in (0,1)$ and for $\gep$ sufficiently small, then
for every $\gb_1\in (0, \gb_0)$ and every positive $\zeta< 2 (1-\gb_1)$ we have 
\begin{equation}
\label{eq:prox}
\lim_{\gep \searrow 0}
\bbP\left( \sup_{t\in [0, \exp(\gep^{-\zeta})]}
\dist(X^\gep_t, M) \le \gep ^{\gb_1} \right) \, =\, 1\, .
\end{equation} 
\end{proposition}
\medskip

Proposition~\ref{th:prox} can be strengthened in a variety of ways. In particular
the next statement
 says that the trajectories naturally stay closer to $M$ than {\sl just} $\gep^{\gb_1}$ (of Proposition~\ref{th:prox}), in fact $\gb_1$ can be chosen larger than $\gb_0$ and even
arbitrarily close to $1$ if we accept to wait for some time so that  the  dynamics drives 
the trajectory closer to $M$.
The rough estimate is that this distance goes 
like $\gep^{\gb_0}\exp(-\gamma t)$, up to the moment in which this distance enters the fluctuation regimes, that is when
it enters a $O(\gep)$ neighbor of the manifold.
It is therefore clear that getting to a distance
$O(\gep)$ from $M$ takes a time $c \vert \log \gep\vert$ for a sufficiently large  $c>0$. 

\medskip
\begin{proposition}
\label{th:contract}
Assume that there exists $\gb_0\in (0,1)$ such that $\dist(X^\gep_0, M) \le \gep ^{\gb_0}$ for $\gep$ sufficiently small.
Then  if $\gb_1\in (\gb_0, 1)$ and $\zeta< 2(1-\gb_1)$
\begin{equation}
\label{eq:contract}
\lim_{\gep \searrow 0}
\bbP\left( \sup_{t\in [c \vert \log \gep \vert, \exp(\gep^{-\zeta})]}
\dist(X^\gep_t, M) \le \gep ^{\gb_1} \right) \, =\, 1\, ,
\end{equation}
where $c:= 2T(\gb_1-\gb_0)/\vert \log(1-\xi_F/2)\vert$.
\end{proposition}

\medskip

\subsection{More on Floquet Theory}
\label{sec:prelim}

 In order to prove Proposition~\ref{th:prox} and  Proposition~\ref{th:contract}
 we  now give more results connected to   Floquet Theory: 
 this part may be viewed as completing  \S~\ref{sec:model}.

\subsubsection{More on Floquet Theory} We restart from \eqref{eq:later}.
The basic idea of the Floquet theory is to consider, for some chosen phase $\theta_0$, a
matrix $Q(\theta_0)$ such that
\begin{equation}
 \Pi(\theta_0+T,\theta_0)\, =\, e^{-T Q(\theta_0)}\, .
\end{equation}
Remark that $Q(\theta_0)$ may not be unique, but always exists since
$\Pi(\theta_0+T,\theta_0)$ is invertible.
Then, by defining $N(\theta_0+t,\theta_0):=\Pi(\theta_0+t,\theta_0)e^{tQ(\theta_0)}$, one obtains the 
following decomposition of the principal matrix solution:
\begin{equation}\label{eq:Floquet decomp 0}
\Pi(\theta_0+t,\theta_0)\, = \, N(\theta_0+t,\theta_0) e^{-t Q(\theta_0)}\,  .
\end{equation}
It is easy to check that $N(\theta_0+t,\theta_0)$ is $T$ periodic, that is $N(\theta_0+t+T,\theta_0)=N(\theta_0+t,\theta_0)$,
and satisfies $N(\theta_0,\theta_0)=I_d$.
By defining $Q(\theta)=\Pi(\theta,\theta_0)Q(\theta_0)\Pi(\theta_0,\theta)$ and
$N(\theta+t,\theta)=\Pi(\theta+t,\theta)e^{tQ(\theta)}$ one can extend
the decomposition \eqref{eq:Floquet decomp 0} to every phase $\theta$:
\begin{equation}\label{eq:Floquet decomp}
\Pi(\theta+t,\theta)\, = \, N(\theta+t,\theta) e^{-t Q(\theta)}\,  .
\end{equation}
From this one easily sees the central fact,  mentioned in \S~\ref{sec:model}
the matrices $Q(\theta)$ are similar: in fact they satisfy for $\theta_1,\theta_2 \in \bbS_T$ 
\begin{equation}\label{eq:link Q}
 Q(\theta_1)\, =\, \Pi(\theta_1,\theta_2)Q(\theta_2)\Pi(\theta_2,\theta_1)\, ,
\end{equation}
and from this we have the notion of  Floquet exponents.

\subsubsection{Floquet Theory and nonlinear evolution}
If $x_t$ is a solution of \eqref{eq:det} starting at time $0$ from a $x_0$ close
to some point $q_\theta$ of $M$, then one can express
the trajectory $y_t$ defined as
$
 y_t = x_t-q_{\theta+t}
$
in terms of the principal matrix solution: a simple calculation (see \cite{cf:Teschl})
shows that
\begin{equation}
\label{eq:delta0}
 y_t\, = \, \Pi(\theta+t,\theta)y_0 + \int_0^t \Pi(\theta+t,\theta+s) \gd_{\theta+s}(x_s) \dd s\, ,
\end{equation}
where the function $\gd$ is defined as
\begin{equation}
\label{eq:delta}
 \gd_\theta (x)\, =\,  F(x)-F(q_\theta)-dF(q_\theta) (x-q_\theta)\, ,
\end{equation}
and satisfies
$\Vert \gd_\theta (x) \Vert  \leq  C \Vert x-q_\theta \Vert^2$
close to $M$ uniformly in $\theta$ for some constant $C$, due to $C^2$ regularity of $F$.

\subsubsection{Floquet Theory: projections and norms}
We will make use of some families of projection. Let us denote by
$P_{\theta}:\bbR^d\rightarrow \bbR^d$
the projection on the tangent space of $M$ at $q_{\theta}$ (generated
by the eigenfunction $F(q_{\theta})$ of $Q(\theta)$) with kernel the
sum of the other characteristic spaces of $Q(\theta)$,
and define $P^\perp_{\theta}= I_d-P_{\theta}$.
From the relation \eqref{eq:link Q} we obtain directly
the following commutation relation:
\begin{equation}\label{eq:conjugaison Pt} 
 P_{\theta_1}\, =\, \Pi(\theta_1,\theta_2)P_{\theta_2}\Pi(\theta_2,\theta_1)\, .
\end{equation}
This relation shows that $P_\theta$ has the same regularity in $\theta$
as $\theta\mapsto \Pi(\theta,\theta_0)$, which is $C^1$, and moreover that there exists a constant $C_P$
such that (uniformly in $\theta$):
\begin{equation}\label{eq:def CP}
 \Vert P_\theta z\Vert\, \leq\, C_P \Vert z\Vert\, \quad \text{and}\quad 
  \Vert P^\perp_\theta z\Vert\, \leq\, C_P \Vert z\Vert\, . 
\end{equation}

Clearly, the hypothesis that the eigenvalues
on the image of $P^\perp_\theta$
have real part strictly negative implies that
there exist $C_\Pi>0$ and $ \gl>0$ such that
for all $z\in \bbR^d$ and $t>0$ (uniformly in $\theta$):
\begin{equation}\label{eq:bound normal space}
 \Vert \Pi(\theta+t,\theta) P^\perp_\theta z\Vert \, =\, \Vert N(\theta+t,\theta)e^{-tQ(\theta)}P^\perp_\theta z \Vert
 \, \leq\, C_\Pi e^{-\gl t}\Vert P^\perp_\theta z\Vert\, .
\end{equation}
On the other hand, since the tangent space of $M$ at $q_\theta$ is the eigenspace associated to the eigenvalue $0$ of $Q(\theta)$, we simply have for all $z\in \bbR^d$ ans $t>0$
(again, uniformly in $\theta$)
\begin{equation}
 \Vert \Pi(\theta+t,\theta)  z\Vert \,  \leq\, C_\Pi\Vert  z\Vert\, .
\end{equation}
It will be useful in the rest of the paper to consider
a family of norms $\Vert\cdot\Vert_\theta$ such that
there exists $\xi_F>0$ (the subscript $F$ is there to stress that it just depends on $F(\cdot)$) such that
\begin{equation}
\label{eq:norm theta contraction}
 \Vert \Pi(\theta+T,\theta) P^\perp_\theta z\Vert_\theta\, \leq\, \left(1-\xi_F\right)\Vert P^\perp_\theta z\Vert_\theta\, ,
\end{equation}
which means that for the topology induced by the norm $\Vert \cdot\Vert_\theta$ the mapping
$\Pi(\theta+T,\theta)$ is a contraction on the normal space at $q_\theta$.
This can be done in a standard way 
 by considering, for a phase $\theta_0$ and
a decomposition in $\mathcal{M}_d(\bbC)$ in Jordan blocs of $\Pi(\theta_0+T,\theta_0)$, i.e.
\begin{equation}
 \Pi(\theta_0+T,\theta_0)\, =\, U^{-1} J U\, , \quad \text{with}\quad 
 J\, =\,
 \left(
 \begin{array}{cccc}
  J_1 & & & \\
    & J_2 & &\\
    & & \ddots &\\
    & & & J_r
 \end{array}
\right)\,  ,
\end{equation}
the norm $\Vert \cdot \Vert_{\theta_0}$ defined as
\begin{equation}
\label{eq:thetanorm}
 \Vert z \Vert_{\theta_0}\, :=\, \Vert E^\gd U z\Vert\, ,
\end{equation}
where the matrix $E^\gd$ has the same structure as $J$ but with the blocs $J_i$
replaced by the blocs of the same size
$E^\gd_i=\text{diag}(\gd^{i-1},\gd^{i-2},\ldots\, ,\gd, 1)$, and by taking $\gd$ small enough.
Remark that one can choose a family of norm $\Vert\cdot\Vert_\theta$ that is smooth
with respect to $\theta$: since $\Pi(\theta+T,\theta)=(\Pi(\theta_0,\theta))^{-1}\Pi(\theta_0+T,\theta_0)\Pi(\theta_0,\theta)$,
one can define $\Vert z\Vert_{\theta}:=\Vert E_\gd \Pi(\theta_0,\theta)Uz\Vert$ with the same small enough
$\gd$ for all $\theta$.

Of course 
if the matrix $\Pi(\theta_0+T,\theta_0)$ is diagonalizable \eqref{eq:thetanorm}  holds with $E^\gd$ equal to the  identity matrix.

%
%
%
%
%
%
%
%
%
%

\subsection{The proof of Proposition~\ref{th:prox} and Proposition~\ref{th:contract}}
The proof starts by a part that is common to both propositions.

\subsubsection{Iterative scheme}
We are interested in the dynamics of the process $X=X^\gep$ on time intervals
of the type $[0,\exp(\gep^{-\zeta})]$.
We will consider for $1\leq k\leq n_f$, with $n_f:= \left\lfloor \frac{\exp(\gep^{-\zeta})}{T}\right\rfloor$, the projections
\begin{equation}
 \ga_{k-1}\, :=\, \mathrm{proj}_M(X_{(k-1)T})\, ,
\end{equation} 
when $X_{(k-1)T}$ is close enough to $M$ such that the projection $\mathrm{proj}_M(\cdot)$
is well defined. We would like to
compare, on each time interval $[(k-1)T, kT]$, the process $X$ to the periodic solution starting from
$q_{\ga_{k-1}}$. In view of what we want to prove and, accessorily 
the projections to be well defined, we introduce stopping times.
More precisely we fix a $\gd>0$ that guarantees that $\text{proj}_M(\cdot)$ is  well defined 
in $M_\gd$ (cf. \eqref{eq:Mdelta}) and
we consider the stopping time, in fact the stopping couple $(k_\tau, \tau)$, such that
\begin{equation}
 (k_\tau,\tau)\, =\, \inf \left\{ (k,s)\in \{1,\ldots,n_f\}\times 
 [0,T]:\, \Vert X_{(k-1)T+s}-q_{\ga_{k-1}+s}\Vert \geq \gd \right\}\, ,
\end{equation}
where the infimum is taken with respect to the lexicographic order.

Define moreover for $k=1,\ldots,n_f$ the times
\begin{equation}
 \tau^k\, :=
 \left\{ 
 \begin{array}{ll}
 T & \text{if}\quad k<k_\tau\, ,\\
 \tau & \text{if} \quad k=k_\tau\, ,\\
 0 & \text{if} \quad k>k_\tau\, .
 \end{array}
\right.
\end{equation}
We can now consider the process defined by $\theta_0=\mathrm{proj}_M(X_0)$, assume
dist$(X_0,M)\le \gd$, and for $1\leq k\leq n_f$ 
\begin{equation}
 \theta_{k}\, :=\, \mathrm{proj}_M(X_{(k\wedge k_\tau-1)T+ t\wedge \tau^k})\, ,
\end{equation}
so that $\theta_k$ is the projection on $M$ of $X_{kT}$, unless this process has been stopped
and in this case it is its projection at the stopping time.
We consider also the process $Y^k$ defined for $t\in [0,T]$ 
and $1\leq k\leq n_f$ by
\begin{equation}
 Y^k_t\, :=\, X_{(k\wedge k_\tau-1)T+ t\wedge \tau^k}-q_{\theta_{k-1}}\, .
\end{equation}
This process $Y^k_t$ satisfies the following mild equation for $t\in[0,T]$ (recall \eqref{eq:delta0}-\eqref{eq:delta}): 
\begin{multline}\label{eq:Y reduced}
 Y^{k}_t\, =\, \Pi(\theta_{k-1}+t\wedge \tau^k,\theta_{k-1})Y^{k}_0 \\
 +\int_0^{t\wedge \tau^k} \Pi(\theta_{k-1}+t\wedge \tau^k,\theta_{k-1}+s)\gd_{\theta_{k-1}
 +s}(X_{(k-1) T+s})\dd s \\
 +\gep\int_0^{t\wedge \tau^k} \Pi(\theta_{k-1}+t\wedge \tau^k,\theta_{k-1}+s)
 G(X_{(k -1)T+s})\dd B_{(k-1)T+s}\, .
\end{multline}
We set
for $t\in[0,T]$ and $k=1,\ldots,n_f$ 
\begin{equation}\label{def:Zk}
 Z^{k}_{t}\, =\, \int_{0}^{t\wedge \tau^k} \Pi(\theta_{k-1}
 +t\wedge \tau^k,\theta_{k-1}+s) G(X_{(k-1)T+s}) \dd B_{(k-1)T+s}\, .
\end{equation}

\subsubsection{Two lemmas}
The iteration we just introduced is controlled via a bound on  the noise term and 
by exploiting the contracting properties of the dynamics neat $M$. This two ingredients 
correspond to the two lemmas we state and prove next.

\medskip
\begin{lemma}
\label{th:ctrl_noise}
Given the stochastic evolution \eqref{eq:main} (hence, given $F(\cdot)$ and $G(\cdot)$) we set
$C_{F,G}:=  T \sup_{\theta\in \bbS_T}\Vert \Pi(0,\theta)\Vert^4 m \Vert G \Vert_\infty^2$.  
For every $q\ge 0$ and every $n_f$
\begin{equation}
\label{eq:ctrl_noise}
\bbP\left( \sup_{k\in 1, \ldots, n_f} \sup_{t\in [0,T]}
\left \Vert Z_t^k \right\Vert \, \ge \, q\right) \, \le \, 2d \, n_f \exp\left( -\frac{q^2}{2C_{F,G}}\right)\, .  
\end{equation}
\end{lemma}
\medskip

\noindent
{\it Proof.}
For a fixed $\alpha\in \bbS_T$ define for $k=1,\ldots,n_f$, $t\in [0,T]$ 
\begin{equation}
 \tilde Z^{k}_{t}\, =\, \int_{0}^{t\wedge \tau^k} \Pi(\alpha,
 \theta_{k-1}+s) G(X_{(k-1)T+s}) \dd B_{(k-1)T+s}\, ,
\end{equation}
so that $Z^k_t=\Pi(\theta_{k-1}+t\wedge t^\tau,\alpha)\tilde Z^k_t$.
Since $\Vert Z^k_t\Vert\leq C_\Pi\Vert \tilde Z^k_t\Vert$, with $C_\Pi := \sup_{\theta\in \bbS_T}\Vert \Pi(0,\theta)\Vert$, 
the left hand side in \eqref{eq:ctrl_noise} is bounded by 
\begin{equation}
\label{eq:sts5}
 \bbP\left(\sup_{k=1,\ldots,n_f}\sup_{t\in[0,T]}\max_{j=1, \ldots d}\left \vert  \left(\tilde Z^k_t\right)_j
 \right\vert\geq \frac{q}{C_\Pi}\right)\, \le n_f d \max_{j=1, \ldots d} \bbP\left(\sup_{t\in[0,T]}\left \vert  \left(\tilde Z^k_t\right)_j
 \right\vert\geq \frac{q}{C_\Pi}\right)
 \, .
\end{equation} 
Therefore by proceeding precisely like in \eqref{eq:expmartGb1}-\eqref{eq:expmartGb2}: the quadratic variation of the $j^{\textrm{th}}$ component of
 $\tilde Z^{k}_{t}$ equals 
 \begin{equation}
  \int_{0}^{t\wedge \tau^k} \sum_{j'=1}^m\left(\Pi(\alpha,
 \theta_{k-1}+s) G(X_{(k-1)T+s})\right)_{j,j'}^2 \dd s \, \le \, mT C_\Pi^2 \left \Vert G \right \Vert^2_\infty\, ,
 \end{equation}
the proof is complete. 
\qed

\medskip

Recall that $\xi_F\in (0,1)$ is given in \eqref{eq:norm theta contraction}. We recall also that the subscript 
$F$ is to stress that such a constant depends only on $F(\cdot)$ and this convention is used also in the next statement.
In fact the constants may also depend on $\gd$, that is the constant entering the definition of the stopping couple, but $\gd$ is fixed and
depends ultimately only on $F(\cdot)$ so we will not write it explicitly.

\medskip

\begin{lemma}
\label{th:iterlem}
Choose a $\underline{\gb} \in (0,1)$.
There exist $\gep_{F,\underline{\gb}}>0$ and $C_F\ge 1$ such that if
for $k\in \bbN$ we have that for every $\gep \in (0, \gep_{F,\underline{\gb}}]$ and a $\gb\ge \underline{\gb}$
\begin{equation}
\label{eq:iterlem-1}
\sup_{t \in [0, T]} \left \Vert Z^k_t \right\Vert_{\theta_{k-1}}\, \le \, \frac{\xi_F}4 \gep^{-1+ \gb}
\ \ \ \ 
\text{ and } \ \ \ \
\left \Vert Y_0^k \right\Vert_{\theta_{k-1}}\, \le \, \gep^{\gb}\, ,
\end{equation}
then, for the same values of $\gep$ and $\beta$, we have that
\begin{equation}
\label{eq:iterlem-2}
\sup_{t\in [0,T]} \left \Vert Y_t^k \right\Vert_{\theta_{k-1}}\, \le \, C_F \gep^{\gb}
\ \ \ \ 
\text{ and } \ \ \ \
\left \Vert Y_0^{k+1} \right\Vert_{\theta_{k}}\, \le \, \left( 1- \frac{\xi_F}2\right) \gep^{\gb}\, .
\end{equation}
\end{lemma}

\medskip

\noindent
{\it Proof.} We start by observing that the second statement in \eqref{eq:iterlem-1} guarantees that $k_\tau>k$, if $\gep_{F, \underline{\beta}}$ is chosen sufficiently small.
We then introduce the stopping time $\tau:=\inf\{ t\ge 0: 
 \Vert Y_t^k \Vert_{\theta_{k-1}}\, > \, C_F \gep^{\gb}\}$: $C_F$ a positive constant that is going to be chosen just below. 
 Note that $\tau < \tau^k$, at least if $\gep_{F, \underline{k}}$ is chosen sufficiently small.
For all $t\leq \tau$ we obtain, using the mild formulation
\eqref{eq:Y reduced}, the upper bound on the normal space \eqref{eq:bound normal space},
the norm equivalence $c_M \Vert \cdot\Vert_\theta
\leq \Vert\cdot\Vert\leq C_M\Vert\cdot\Vert_\theta$ and the fact that $\gd_\theta$ is quadratic around $M$ (cf. \eqref{eq:delta}:
the positive constant associated to this quadratic bound is denoted by $c_{F,2}$), we obtain
\begin{multline}
\label{eq:estim Vert Yt}
 \Vert Y^k_t\Vert_{\theta_{k-1}} \, \leq\, \Vert \Pi(\theta_{k}+t,\theta_{k})Y^{k}_0\Vert_{\theta_{k-1}}
 +\int_0^t \Vert \Pi(\theta_{k}+t,\theta_{k}+s)\gd_{\theta_{k}
 +s}(X_{k T+s}) \Vert_{\theta_{k-1}}\dd s
 +\gep \Vert Z^k_t\Vert_{\theta_{k-1}}\\
 \, \leq\,c_M C_M  C_\Pi \big( e^{-\gl t} \gep^{\gb } +  T c_{F,2} C_F^2 \gep^{2\gb}\big)+\frac{\xi_F}4 \gep^{\gb}\,
 \le \, \left(2  c_M C_M C_\Pi + \frac{\xi_F}4\right) \gep^{\gb}\, ,
\end{multline}
where in the second inequality we have chosen $\gep_{F,\underline{\gb}}$ so that $T c_{F,2} C_F^2  (\gep_{F,\underline{\gb}})^{\underline{\gb}}\le 1$.
Note that, with the choice $C_F= 2  c_M C_M C_\Pi+ \xi_F/4$ and possibly by choosing $\gep_{F,\underline{\gb}}$ even smaller, \eqref{eq:estim Vert Yt} implies $\tau \ge T$,
which proves \eqref{eq:iterlem-2}.

It remains now to show that
$\Vert Y^{k+1}_0\Vert_{\theta_k}\leq \gep^{\gb}$. We restart from 
\eqref{eq:Y reduced}, apply to all terms the projection $P^\perp_{\theta_{k-1}}$ and then 
again apply the norm, getting to the analog  of the first line of \eqref{eq:estim Vert Yt}
 for $\Vert P^\perp_{\theta_{k-1}} Y^k_{T}\Vert_{\theta_{k-1}}$.
Using \eqref{eq:norm theta contraction} we obtain, in strict analogy with the second step in \eqref{eq:estim Vert Yt}
(and choosing $\gep_{F, \underline{\beta}}$ small enough), that
\begin{equation}
\label{eq:Port}
\Vert P^\perp_{\theta_{k-1}} Y^k_{T}\Vert_{\theta_{k-1}} \, \leq\,   \left(1-\xi_F\right)\gep^{\gb} 
+ c_{F,3}\gep^{2\gb}+\frac{\xi_F}4\gep^{\gb}\, \le\,  \left(1-\frac 34 \xi\right)\gep^{\gb }\,.
\end{equation}
To deduce an upper bound for $Y^{k+1}_0$ with this result remark that one may write
$Y^{k+1}_0$ in terms of $P^\perp_{\theta_{k-1}}Y^k_{T}$:
\begin{multline}
\label{eq:decomp Ykplus10}
 Y^{k+1}_0\, =\, q_{\theta_{k-1}} + Y^{k}_{T}-q_{\theta_k}
 \, =\, P^\perp_{\theta_k}\big( q_{\theta_{k-1}} + Y^{k}_{T}-q_{\theta_k}\big)\\
 =\,P^\perp_{\theta_k}\left( q_{\theta_{k-1}} -q_{\theta_k}\right)
 +\big(P^\perp_{\theta_{k}}-P^\perp_{\theta_{k-1}}\big)Y^{k}_{T}+P^\perp_{\theta_{k-1}}Y^k_{T}\, .
\end{multline}
Observe now that
since $\theta_k-\theta_{k-1}=\mathrm{proj}_M(X_{kT})-\mathrm{proj}_M(X_{(k-1)T})$, 
the smoothness of the projection $\mathrm{proj}_M(\cdot)$ implies that
$|\theta_k-\theta_{k-1}|\leq c_{F,4}\Vert X_{kT}-X_{(k-1)T}\Vert \leq c_{F,4}(\Vert Y^k_0\Vert +\Vert Y^k_{T}\Vert)\le
c_{F,4}(C_M + C_F) \gep^\gb$. 
By Taylor expansion
\begin{equation}
P^\perp_{\theta_k}\left( q_{\theta_{k-1}} -q_{\theta_k}\right)\, =\, (\theta_{k-1}-\theta_k)P^\perp_{\theta_k}q'_{\theta_k}
+O\big((\theta_k-\theta_{k-1})^2\big)\, ,
\end{equation}
where once again the term $O\big((\theta_k-\theta_{k-1})^2\big)$ is bounded by $(\theta_k-\theta_{k-1})^2$ times
a constant that depends only on $F(\cdot)$. Since $P^\perp_{\theta_k}q'_{\theta_k}=0$ we conclude
that 
\begin{equation}
\left\Vert P^\perp_{\theta_k}\left( q_{\theta_{k-1}} -q_{\theta_k}\right) \right\Vert_{\theta_{k-1}}
\, \le \, c_{F,5}\gep^{2\gb}\, .
\end{equation}
A very similar estimate holds for  the second term of the right-hand side of \eqref{eq:decomp Ykplus10}, since the smoothness of the application
$\theta \mapsto P^\perp_\theta$ implies that
\begin{equation}
 \big\Vert \big(P^\perp_{\theta_{k}}-P^\perp_{\theta_{k-1}}\big)Y^{k}_{T}\big\Vert_{\theta_{k-1}}
 \, \leq\, c_{F,6}|\theta_k-\theta_{k-1}|\big\Vert Y^{k}_{T}\big\Vert\, \le\, c_{F,7}\gep^{2\gb}\, .
\end{equation}
So one has $\Vert Y^{k+1}_0\Vert_{\theta_{k-1}} \leq (1-2\xi/3)\gep^{\gb}$ when $\gep $ smaller than a constant
that depends only on $F(\cdot)$ and $\underline{\gb}$
and to conclude the proof remark that
\begin{multline}
 \Vert Y^{k+1}_0\Vert_{\theta_{k}}\, =\,\Vert Y^{k+1}_0\Vert_{\theta_{k-1}} +
 \left(
 \Vert Y^{k+1}_0\Vert_{\theta_{k}}-\Vert Y^{k+1}_0\Vert_{\theta_{k-1}}
 \right)
  \\
 \le \,  
 \Vert Y^{k+1}_0\Vert_{\theta_{k-1}}+c_{F,8}\left|\theta_k-\theta_{k-1}\right|\left\Vert Y^{k}_{T}\right\Vert
\,
  \le \, \Vert Y^{k+1}_0\Vert_{\theta_{k-1}}+c_{F,9}\gep^{2\gb}\\ \le \, \left( 1- \frac{\xi_F}2\right)\gep^{\gb}\, ,
\end{multline}
since the mapping $\theta\mapsto \Vert\cdot\Vert_\theta$ is smooth and, again, we have chosen $\gep$ suitably small. This 
establishes \eqref{eq:iterlem-2} and completes the proof of Lemma~\ref{th:iterlem}.
\qed

\medskip

\noindent
{\it Proof of Proposition~\ref{th:prox}.}
Since $\gb_1< \gb_0$ and $\zeta < 2(1-\gb_1)$, we can choose $\gb_2 \in (\gb_1, \gb_0)$ such that
$\zeta < 2(1-\gb_2)$ still holds. 
 The second inequality in \eqref{eq:iterlem-1} holds
for $\gb= \gb_2$ and  $\gep$ sufficiently small. The first one holds, again for $\gb= \gb_2$
and, by Lemma~\ref{th:ctrl_noise}, for every $n=1,2 , \ldots, \lfloor \exp(\gep^{-\zeta})/T \rfloor$ 
with probability going to one as $\gep$ tends to zero. 
The conclusion is of course \eqref{eq:iterlem-2} for $k=1$. It is then clear that this procedure can be iterated
up to $k= \lfloor \exp(\gep^{-\zeta}) /T\rfloor$ and leads to 
\begin{equation}
\lim_{\gep \searrow 0}
\bbP\left( \sup_{t \in [0,  T(1+\lfloor \exp(\gep^{-\zeta}) /T\rfloor]}
\left \Vert Y_t^k \right\Vert_{\theta_{k-1}} \le \, C_F 
\gep^{\gb_2} \right)\, =\, 0\,,
\end{equation}
which implies
\begin{equation}
\lim_{\gep \searrow 0}
\bbP\left( \sup_{t \in [0,   \exp(\gep^{-\zeta}) ]}
\text{dist} (X_t^\gep,M) \,\le \, C_M C_F 
\gep^{\gb_2} \right)\, =\, 0\,,
\end{equation}
and \eqref{eq:prox} is proven because $\gb_1$ is smaller than $\gb_2$. 
\qed 

\medskip

\noindent
{\it Proof of Proposition~\ref{th:contract}.}
Choose $\gb_0^-\in (0, \gb_1)$ and $\gb_1^+\in (\gb_1,1)$ in such a way that $\zeta$, which is smaller than
$2(1-\gb_1)$, is also smaller than $2(1-\gb_1^+)$ and $\gb_1^+-\gb_0^- < 2( \gb_1-\gb_0)$.
In particular, with these choices we have that $\dist(X^\gep_0, M) \le \gep ^{\gb_0}$
implies $\Vert Y^1_0 \Vert_{\theta_0} \le \gep^{\beta_0^-}$ for $\gep$ small.
By Lemma~\ref{th:ctrl_noise} we directly obtain that
\begin{equation}
\label{eq:from-ctrl_noise}
\bbP\left( \sup_{k\in 1, \ldots, \lceil \exp( \gep^{-\zeta})/T\rceil} \sup_{t\in [0,T]}
\left \Vert Z_t^k \right\Vert \, \ge \, \frac{\xi_F}4\gep^{-1+ \gb_1^+}\right) \, = \,0\, .  
\end{equation}
Let us set $\eta_F:= 1-\xi_F/2 \in (0,1)$. 
For  $k=0,1, \ldots$
 we introduce 
\begin{equation}
\gb_k\, =\, \gb^-_0+k \frac{\log \eta_F}{\log \gep}\, ,
\end{equation}
so $\gep^{\gb_k}:= \gep ^{\gb_0^-}
\, \eta_F^k$, and we set $k_\gep:= \max\{k:\, \gep^{\gb_k}\ge \gep^{\gb_1^+}\}$ so
\begin{equation}
\label{eq:k-eps}
k_\gep\, =\, \left\lfloor \left(\gb_1^+ -\gb_0^-\right)\frac{\log \gep}{\log \eta_F}\right\rfloor\, \le \,
 2\left(\gb_1 -\gb_0\right)\frac{\log \gep}{\log \eta_F}\, ,
\end{equation}
where the inequality holds for $\gep$ small.
Then we observe that, by \eqref{eq:from-ctrl_noise}, the first condition in \eqref{eq:iterlem-1} is satisfied
for $k=0,1, \ldots, k_\gep$. Moreover by definition of $\gb_k$
\begin{equation}
\left \Vert Y_0^{k+1} \right\Vert_{\theta_{k}}\, \le \, \left( 1- \frac{\xi_F}2\right) \gep^{\gb_k} \ \ \Leftrightarrow \ \
\left \Vert Y_0^{k+1} \right\Vert_{\theta_{k}}\, \le \, \gep^{\gb_{k+1}}\, ,
\end{equation}
so  the second inequality in \eqref{eq:iterlem-2} with $\gb=\gb_k$ is 
is the second inequality in \eqref{eq:iterlem-1} with $k$ replaced by $k+1$ and $\gb=\gb_{k+1}$.
So Lemma~\ref{th:iterlem} can be iterated $k_\gep+1$ times and
we retain that $\Vert Y ^{k_\gep+1} \Vert_{\theta_{k_\gep}} \le \gep^{\gb_1^+}$. 

We now can apply 
Lemma~\ref{th:iterlem}, $\lceil \exp( \gep^{-\zeta})/T\rceil -k_\gep$ times without taking advantage 
of the contraction factor $\eta_F$ (our noise estimate \eqref{eq:from-ctrl_noise} is not good enough to
get an advantage from this contraction factor), so we just keep the estimate $\Vert Y ^{k_\gep+1} \Vert_{\theta_{k_\gep}} \le \gep^{\gb_1^+}$
up to $k=\lceil \exp( \gep^{-\zeta})/T\rceil +1$, as well as the first inequality in Lemma~\ref{eq:iterlem-2}. And it is
precisely the first inequality in Lemma~\ref{eq:iterlem-2} that yields -- once we  recall 
\eqref{eq:k-eps}  -- \eqref{eq:contract}. The proof of
Proposition~\ref{th:contract} is complete.
\qed

\section{Completion of the proof of Theorem~\ref{th:main2}}
\label{sec:complete}

In view of what we want to prove let us choose $t_\gep\in[ a \exp(\gep^{-\zeta}), \exp(c \gep^{-2})]$, for some $a>0$,
and we set $\tau_f:= \tau_{M_\gd^\complement}$ to make formulas more compact.
We now  consider again the processes $Z^k$ defined in \eqref{def:Zk} and 
we use the notation
\begin{equation}
\label{eq:Qk}
Q^k_t \, :=\, \int_{(k-1)T}^{(k-1)T+t} G\left(X^\gep_{s \wedge \tau_f}\right) \dd B_s \, . 
\end{equation}
We then proceed  by defining the following family of  integer valued stopping times (the filtration $\{\cF_t \}_{t\ge 0}$
 is always the Brownian one): $k^0=0$ and for $i\geq 1$
\begin{equation}
k^i\,=\, \inf\left\{ k> k^{i-1}: \, \sup_{t\in [0,T]}\left\Vert Z^{k}_t\right\Vert 
\vee \left\Vert Q^{k}_t\right\Vert 
\, \ge\,  \gep^{\beta-1}
\right\}\, ,
\end{equation}
where $\gb$ is chosen so that $2(1-\gb)< \zeta$, so if $\xi < 2(1-\gb)$ then $\exp(\gep^{-\xi}) \ll t_\gep$.
These stopping times correspond to the {\sl rather} large excursions experienced either by  $Z^k_\cdot$
or $Q^k_\cdot$, which can induce large excursions for $X^\gep_\cdot$. These large excursions should be very rare, and be most of the time separated by a very large time interval, so that the process $X^\gep$ has plenty of time to relax and stay in a neighborhood of $M$ before the next excursion.
More precisely, for $\xi\in (0,2(1-\beta))$, let us define the following random variable, denoting
\begin{equation}
D^1_\gep\, :=\, \#\left\{i\in \bbN:\,  k^i\leq n_f \text{ and } k^i-k^{i-1}\leq e^{\gep^{-\xi}}\right\}\ \ \ \text{ with } \ \ \
n_f\,:=\,\lfloor t_\gep/T\rfloor\, .
\end{equation}
We also introduce $D^0_\gep\, :=\, \#\{i\in \bbN:\,  k^i\le n_f\}$.
\medskip

\begin{lemma}
\label{lem:bound Dgep}
For all $\xi'\in (\xi,2(1-\beta))$ we have
\begin{equation}
\label{eq:bound Dgep}
\lim_{\gep\searrow 0}\bbP\left(D_\gep^0 < t_\gep e^{-\gep^{-\xi'}}\text{ and }  D^1_\gep < t_\gep e^{\gep^{-\xi}-2\gep^{-\xi'}} \right)\, =\, 1\, .
\end{equation}
\end{lemma}

\medskip

\noindent
{\it Proof.} Recall that both $Z^k_t$ and $Q^k_t$ are measurable with respect to $\cF_{(k-1)T+t}$.
Lemma \ref{th:ctrl_noise}  implies that for $\xi''\in (\xi,2(\beta-1))$
\begin{equation}
\bbP\left(\sup_{t\in[0,T]} \left\Vert Z^k_t\right\Vert\vee \left\Vert Q^k_t\right\Vert> \gep^{\beta-1}\bigg \vert \cF_{(k-1)T}\right)\, \leq\,  e^{-\gep^{-\xi''}}\, =:\, p_\gep\, ,
\end{equation}
so the point process $\{k^i\}_{i= 0, 1, \ldots }$ is stochastically dominated by a renewal process $\{\kappa^i\}_{i= 0, 1, \ldots }$ 
with geometric inter-arrival of parameter $p_\gep$. Therefore $\bbE (D_\gep^0 ) \le n_f p_\gep$ and the probability 
that $D_\gep^0$ is larger than $t_\gep e^{-\gep^{-\xi'}}$ vanishes by Markov inequality. 
Define now the process $\{\cN_j\}_{j= 0, 1, \ldots }$ by $\cN_0:=0$ and 
\begin{equation}
\cN_j\, :=\, \sum_{i=1}^j \ind_{\kappa^i-\kappa^{i-1} \leq e^{\gep^{-\xi}}}\, ,
\end{equation}
and the stopping time
\begin{equation}
\iota \, :=\, \inf\{j=0,1, \ldots:\, \cN_j>n_f\}\, .
\end{equation}
Denote $q_\gep:=\bbP(\kappa^1\leq e^{\gep^{-\xi}})$. Since $\lim _{\gep \searrow 0}
p_\gep \exp(\gep^{-\xi})=0$ we have $q_\gep\sim p_\gep e^{\gep^{-\xi}} =e^{\gep^{-\xi}-\gep^{-\xi''}}$, and the process 
$\{\cN_j- q_\gep j\}_{j=0,1, \ldots}$ is a martingale, so $\bbE (\cN_\iota)= q_\gep\bbE (\iota)$. Remark moreover that $\iota-1$ simply follows a binomial distribution of parameters $n_f$ and $p_\gep$, so that $\bbE(\iota)=p_\gep n_f +1$. Finally 
$\bbE (\cN_\iota)\sim n_f p_\gep q_\gep $, and we conclude by Markov inequality, once the choice of  a $\xi'\in (\xi,\xi'')$
is made.
\qed

\medskip

\noindent
{\it Proof of Theorem~\ref{th:main2} assuming \eqref{eq:cond-main2.2}.}
%
Let us observe that in Lemma~\ref{lem:bound Dgep} the upper bound on $t_\gep$, namely the constant $c>0$, 
is arbitrary. But the lemma carries little information 
if $X^\gep_\cdot$ exits $M_\gd$, because both the $Z$ and $Q$ processes depend on $X^\gep_\cdot$ stopped at $\tau_f$. 
However
by a standard Large Deviations argument (based on  \cite[Th.~5.6.7]{cf:DZ}) one establishes that there exists $c>0$
such that 
\begin{equation}
\lim_{\gep \searrow 0} \bbP \left( \tau_f \ge \exp\left( c \gep^{-2} \right)\right) \, =\, 1\, .
\end{equation}
This identifies the constant $c>0$ appearing in Theorem~\ref{th:main2} and one can make it more explicit in terms of  the so called 
 {\sl quasi-potential}  associated to exiting $M_\gd$, cf. \cite{cf:FW}. We therefore choose to work on the event $\tau_f \ge \exp\left( c \gep^{-2} \right)$. 
 
 It is also practical to exploit another (rough) Large Deviations estimate, namely that there exists
 $A>0$
  \begin{equation}
  \label{eq:roughLD2}
 \lim_{\gep \searrow 0}
\bbP\left(\sup_{k=1,\ldots,n_f}\sup_{t\in[0,T]} \left\Vert Q^k_t\right\Vert\, \ge \,  A\gep^{-1}\right)\, =\, 0\, .
\end{equation}
 A proof of this claim is directly achieved by 
 applying  Lemma \ref{th:ctrl_noise} and, given $G(\cdot)$  and $M_\gd$, how large $A$ has to be chosen only depends on the value of $c$: just use that in any case $n_f \le \exp(c\gep^{-2})/T$.

 Therefore we can and will work assuming (1) that the event whose probability is estimated in \eqref{eq:bound Dgep} is realized, (2) that 
$\tau_f \ge \exp\left( c \gep^{-2} \right)$ and (3) that 
 the complement of the event whose probability is evaluated in
 \eqref{eq:roughLD2} is realized. There are three cases to treat:
\smallskip

\begin{enumerate}
\item The bad blocks $[(k^i-1)T, k^iT]$. 
\item The {\sl short gaps} between bad blocks, that is the intervals $[k^{i-1}T, (k^i-1 )T]$ that are not longer
than $Te^{\gep^{-\xi}}$.
\item The {\sl long gaps} between bad blocks, i.e. what is left.
\end{enumerate} 

\smallskip

The first two cases are treated simply by exploiting the rough estimate \eqref{eq:roughLD2} which warrants 
that on time intervals of length $t$ the phase cannot change more than $Ct$, with $C$ a constant that depends on 
$A$ and on the drift $F(\cdot)$: $F(\cdot)$ is bounded because we are bound to $M_\gd$, since $\tau_f \ge \exp\left( c \gep^{-2} \right)$. Therefore, recalling Lemma \ref{lem:bound Dgep}, the total contribution of the bad blocks to the phase does not exceeds $Ct_\gep \exp(-\gep^{-\xi'})$ and 
the contribution of the short gaps does not exceed $Ct_\gep \exp(2(\gep^{-\xi}-\gep^{-\xi'}))$. Since these two quantities are $o(t_\gep)$, they give no contribution to the 
final result.

The leading contribution comes all from the long gaps.
For this part we use the sharper estimates on the noise given by the second requirement 
in the event whose probability is evaluated in \eqref{eq:bound Dgep}. This tells us that the noise term is small
in a way that we can adapt  the proof of Theorem~\ref{th:main} -- much like we did for the proof of Theorem~\ref{th:main2} assuming \eqref{eq:cond-main2.1} at the end of Section~\ref{sec:proofs} -- and of
Proposition~\ref{th:prox}. Let us see this step more in details. 

First remark that for each of these long gaps the initial point $X^\gep_{k^{i-1}T}$ may be anywere in $M_\gd$, and we have to show that the process comes back quickly in $M_{\beta'}$ for $\beta'<\beta$. This can be done exactly as the step 1 of proof of Theorem \ref{th:main} and in the proof of Proposition \ref{th:contract}, showing that it happens in a time smaller than $-C\log\gep$ with $C$ depending only on $\gd$ and $\beta'$.
So during this period of time the phase cannot change more than $-C'\log\gep=o((k^i-k^{i-1})T)$.
This first part can thus be neglected, and we can place ourself in the case $X^\gep\in M_{\gep^{\beta'}}$. We follow then what has already been done at the end of Section \ref{sec:proofs}. We want to prove that
\begin{equation}
\tilde{\theta \left(X^\gep_{\cdot\wedge \tau_f}\right)}\left( (k^i-1) T \right) -
\tilde{\theta \left(X^\gep_{\cdot\wedge \tau_f}\right)}\left( k^{i-1} T \right) 
-t_k\\ =\, b \gep^2 t_k +o\left(\gep^2 t_k\right)\, ,
\end{equation}
with $t_k=\left(k^i-k^{i-1}\right)T$, and we thus consider, similarly as in \eqref{eq:ut again},
\begin{equation}
u^\gep(t) \, =\, 
\tilde{\theta \left(X^\gep_{\cdot\wedge \tau_\gep}\right)}\left( t \right) -
\tilde{\theta \left(X^\gep_{\cdot\wedge \tau_f}\right)}\left( k^{i-1} T \right)-
\tilde{\theta \left( x_{\cdot\wedge \tau_\gep}\right)}\left( t\right)\, ,
\end{equation}
with $x_t=\Phi(X^\gep_{k^{i-1}T},t)$. This can be treated exactly as in the end of Section \ref{sec:proofs}, recalling that the noise (i.e. the $Q^k$'s and the $Z^k$'s) is controlled on this period of time. More precisely the quadratic variation 
term can be neglected like in the partial proof of Theorem~\ref{th:main2} at the end of Section \ref{sec:proofs} 
and the other terms -- steps form \eqref{eq:ito1} to \eqref{eq:fr18} -- require the control of the noise term 
provided precisely by the control on the $Q^k$'s and the fact that $X^\gep_\cdot$ stays in a $\gep^{\gb}$ neighborhood
of $M$, that follows from the control on the $Z^k$'s by applying  Lemma~\ref{th:iterlem}: the details are the same as in the proof of Proposition~\ref{th:prox}. This completes the proof of Theorem~\ref{th:main2}.  
\qed

\section*{Acknowledgements}
This work stems from interactions with various colleagues, among them we mention and thank
Mathieu Merle and Khashayar Pakdaman.
C.P. acknowledges the support of ERC grant MALADY.

\end{document}